\documentclass[12pt, leqno, twoside]{article}
\usepackage{a4, indentfirst, latexsym,amssymb}
\usepackage{t1enc}

  \newtheorem{Th}{Theorem}[section]
  
  \newtheorem{Prop}[Th]{Proposition}

  \newtheorem{Lem}[Th]{Lemma}

  \newtheorem{Rem}[Th]{Remark}
  
  \newtheorem{Ex}[Th]{Example}

  \def\conv{\mathrm{conv}\,}
  \def\d{\,\mathrm{d}}

  \def\deg{\mathrm{deg}}
  \def\Deg{\mathrm{Deg}}
  \def\la{\langle}
  \def\ra{\rangle}

  \def\cl{\o }

  \def\o#1{\overline{#1}}

  \def\part{\partial}

  \def\R{\mathbb{R}}
  
  \def\Q{\mathbb{Q}}

  \newcommand{\lma}{\lambda}
  
  \newcommand{\eps}{\varepsilon}
  \newcommand{\be}{\begin{equation}}
  \newcommand{\ee}{\end{equation}}

\begin{document}
\begin{center}
{\large{\bf Krasnosel'skii type formula and translation along \\ trajectories method
for evolution equations }}

Aleksander \'{C}wiszewski
\footnote{Corresponding author, e-mail: {\tt aleks@mat.umk.pl}

},
Piotr Kokocki\\
{\em Faculty of Mathematics and Computer Science\\
Nicolaus Copernicus University\\
ul. Chopina 12/18, 87-100 Toru\'n, Poland }\\

\end{center}

\begin{abstract}
The Krasnosel'skii type degree formula for the equation $\dot u = - Au  + F(u)$
where $A:D(A)\to E$ is a linear operator on a separable Banach space $E$ such that
$-A$ is a generator of a $C_0$ semigroup of bounded linear operators of $E$ and $F:E\to E$ is a locally Lipschitz $k$-set contraction, is provided. Precisely, it is shown that if $V$ is
an open bounded subset of $E$ such that $0\not\in (-A+F)(\part V \cap D(A))$, then
the topological degree of $-A+F$ with respect to $V$ is equal to the fixed point index of
the operator of translation along trajectories for sufficiently small positive time.
The obtained degree formula is crucial for the method of translation along trajectories.
It is applied to the non-autonomous periodic problem and an average principle is derived.
As an application a first order system of partial differential equations is considered.
\end{abstract}

\section{Introduction}

We are concerned with the periodic problem associated with a differential equation of the form
$$\left\{ \begin{array}{l}
\dot u(t) = - A u(t) +F (t,u(t)) \mbox { on } [0,T]\\
u(0)=u(T) \end{array} \right. \leqno{(P_{A,F})} $$
where $A:D(A)\to E$ is a linear operator on a separable Banach space such that
$-A$ generates a $C_0$ semigroup $\{S_A(t):E\to E\}_{t\geq 0}$ of bounded linear operators such that
$\|S_A(t)\| \leq e^{-\omega t}$, for some fixed $\omega>0$ and
any $t\geq 0$, and $F:[0,T]\times E\to E$ is a locally Lipschitz map with respect to the second variable, having sublinear growth and with
$k\in [0,\omega)$ such that, for any bounded $\Omega\subset E$,
$\beta(F(\Omega)) \leq k \beta(\Omega)$ where $\beta$ stands for the Hausdorff measure of noncompactness.\\
\indent In general, one may speak of two topological approaches to the $T$-periodic problem associated with $(P_{A,F})$ ($T>0$ is a fixed period).
One is to formulate it as a fixed point problem for a proper solution operator in a suitable space of functions (see e.g. \cite{ObuKaZ}) and the second way for finding $T$-periodic solutions, which we follow,
is to seek fixed points of the translation along trajectory operator $\Phi_T:E\to E$ associated with the equation $(P_{A,F})$ (see e.g. \cite{Bothe} and \cite{BaKr2}). The crucial point in our approach is the proper version of  {\em a degree formula} stating that, if $F$ is time-independent, then the fixed point index of the operator $\Phi_t:E\to E$, for small $t>0$, is equal to the topological degree of $-A+F$
(with respect to proper open bounded subsets of $E$).  We shall prove such a degree formula, which is an infinite dimensional extension of the
Krasnosel'skii theorem for ordinary differential equations in finite dimensional spaces (see \cite{Krasnoselskii}).
By use of the obtained formula, we prove a criterion for finding periodic solutions of $(P_{A,F})$ in the spirit of Mawhin (\cite{Mawhin69}), Schiaffino-Schmitt (\cite{Schi-Schm})
and Kamenskii-Obukhovskii-Zecca (\cite{{ObuKaZ}}).
The obtained abstract method is applied to a system of partial differential equations. The paper is organized as follows. In Section 2 we provide
basic definitions and facts concerning the topological degrees for $k$-set contractions
and perturbations of operators generating $C_0$ semigroups of contractions.
In Section 3, we prove the degree formula, which is to be used in one of the next sections,
for the the equation $\dot u = - u +F(u)$, where $F$ is a $k$-set contraction.
Section 4 is devoted to general compactness properties of solution operators
for $\dot u = - Au+F(t,u)$ with initial conditions.
In Section 5 we prove the main result
of the paper -- the degree formula. We reduce it, by a proper homotopy, to the result of
Section 3. Finally, Section 6 provides an abstract continuation principle
for periodic solutions and its application.

\section{Topological degree}

In the first part of this section we provide basic information on the degree theory for
$k$-set contraction, mainly due to \cite{Nussbaum}. In the second part of the section we carry out a standard construction of the degree for perturbations of generators of $C_0$ semigroups of contractions. For a survey on the topological degree theory for perturbations of accretive operators we refer to \cite{Kartsatos}.

Recall that a continuous map $F:X\to E$, where $X$ is a closed bounded subset of a Banach space $E$, is called
a {\em $k$-set contraction} if there exists a constant $k\in [0,1)$ such that, for any
$\Omega\subset X$,
$$
\beta (F(\Omega))\leq k \beta (\Omega)
$$
where $\beta$ stands for the Hausdorff measure of noncompactness
given by $\beta(\Omega): =\inf\{r>0\,|\, \Omega \mbox{ can be covered with a finite number of balls (in $E$) of radius } r \}$
(see e.g. \cite{Akhmerov-etal}, \cite{Deimling} or \cite{ObuKaZ}).
A compact convex set $C\subset E$ is said to be  {\em a fundamental set}
(for $F$) provided $F( X \cap C ) \subset C$ and, for any $x\in X$, the condition
$x \in \o{\conv} ( \{F(x)\} \cup C )$ implies $x\in C$.
\begin{Rem}\label{construct-fundset}{\em
Let $F:X\to E$ be as above and define $( C_n )_{n\geq 0}$  a sequence by
$C_{0}:=\o{\conv} ( \{x_{0}\} \cup F(X) )$,
$C_{n}:= \o{\conv} ( \{x_{0}\} \cup F(X \cap C_{n-1})), \mbox{ for } n\geq 1$
where $x_0$ is an arbitrary point from $X$.
Then, one can verify that, $C_n\subset C_{n-1}$, for $n\geq 1$, and
$C:=\bigcap_{n=1}^{\infty}C_{n}$ is a fundamental set for $F$.}
\end{Rem}
A vector field $I-F:\o U\rightarrow E$, where $U\subset
E$ is an open bounded set, is called {\em admissible} if $F$ is a
$k$-set contraction  (with some $k\in [0,1)$) and $F$ has no fixed points
in the boundary $\part U$.
Let $C$ be any fundamental set for $F$ and $F_C:\o U \to E$ be
a continuous extension of $F_{|\o U\cap C}$
such that $F_C (\o U)\subset \o \conv F(\o U\cap C)$, existing due to the
Dugundji extension theorem (see \cite{Dugundji-Granas}).
Then, by the invariance of $C$, $F_C(\o U)\subset C$, which implies that
$F_C$ is compact. Moreover, if, for some $x\in \o U$, $F_C (x) = x$,
then $x\in C$, which gives
$F(x)=F_C(x)=x$ and $x\in U$. Therefore the following definition of
the topological degree makes sense
\be\label{def-of-deg} \deg ( I -
F , U ) :=  \deg_{LS} (I-F_C,U)
\ee where $\deg_{LS}$ stands for the Leray-Schauder degree (with respect to $0\in E$).
One may prove that the degree is independent of the choice
of a fundamental set and an extension $F_C$.
It has the following standard properties of a topological degree:\\
\noindent {\bf (D1)} If $\deg(I-F,U)\neq 0$, then there exists $x\in U$ such that $F(x)=x$.\\
\noindent {\bf (D2)} If $F:\o U\to E$ is admissible and $U_1, U_2\subset U$ are open and disjoint sets
such that $\{ x\in \o U\,|\, F(x)=x\} \subset U_1 \cup U_2$, then
$$
\deg(I-F, U) = \deg(I-F, U_1)+\deg (I-F,U_2).
$$
\noindent {\bf (D3)} Let $H:\o U\times [0,1]\to E$ be an admissible homotopy, i.e.
a continuous map such that\\
\indent (i) there exists $k\in [0,1)$ such that, for any $\Omega\subset \o U$,
$$
\beta (H(\Omega\times [0,1])) \leq k \beta (\Omega),
$$
\indent (ii)  $H(x,\lma)\neq x$ for any $(x,\lma)\in \part U\times [0,1]$.\\
Then
$$
\deg(I-H(\cdot,0), U)=\deg(I-H(\cdot,1),U).
$$
\noindent {\bf (D4)} If $x_0\in U$, then $\deg(I-x_0, U) = 1.$\\

\begin{Rem} \label{homotopy-compact}
{\em It will be useful in the sequel to know that
if $I-F:\o U\to E$ is an admissible $k$-set contraction vector field, then there exists
a locally Lipschitz compact vector field
$I-F_L:\o U\to E $ homotopic to $I-F$ via an admissible homotopy.
Indeed, find $F_C$ as in (\ref{def-of-deg}).
Since $F$ is a $k$-set contraction, its fixed point set is compact and one has
$\rho:=\inf\{ \|x-F(x)\|\,|\, x\in \part U \}>0.$
By the Lasota-Yorke theorem there exists a locally Lipschitz map $F_L:\o U\to E$ such that
$\|F_L(x)-F_C(x) \| \leq \rho/2$, for any $x\in \o U$ and $F_L(\o U)\subset \o \conv
F_C(\o U)\subset C$. To see that $I-F_L$ is homotopic to $I-F$, define
$H:\o U\times [0,1]\to E$ by
$$
H(x,\lma):= (1-\lma) F(x)+\lma F_L(x).
$$
Clearly, for any  $\Omega \subset \o U$, one gets
$$
\beta(H(\Omega\times [0,1])) \leq \beta (\conv[F(\Omega)\cup C])
= \beta (F(\Omega))\leq k \beta (\Omega).
$$
Furthermore, if for some $(x,\lma)\in \part U\times [0,1]$, $H(x,\lma)=x$, then
$x\in \conv ( \{F(x)\}\cup C )$, which implies $x\in C$ and, consequently,
$\|x-F(x)\| = \lma \| F_L(x) - F (x)\| =\lma \| F_L(x) - F_C (x)\| \leq \rho/2,$
a contradiction.}
\end{Rem}

The topological degree can be extended onto  the class
of {\em condensing vector} fields. For any {\em admissible} condensing
vector field $I-F:\o U\to E$ with open bounded $U\subset E$, i.e. such that
$F$ is continuous, $\beta(F(\Omega)) <\beta(\Omega)$ for any
$\Omega\subset \o U$ with $\beta(\Omega)>0$, and
$0\not\in (I-F)(\part U)$, there exists $\o\lma\in (0,1)$ such
that
\be\label{no-zeros-on-boundary-cond-2}
0\not\in (I-(1-\lma)F)(\part U) \ \ \mbox{ for any } \lma\in (0,\o \lma].
\ee
Indeed, if we suppose to the contrary, then there are
$\lma_n\to 0^+$ and $(x_n)\subset \part U$ such that $x_n = (1-\lma_n)
F(x_n)$, for each $n\geq 1$. Then, either $\beta(\{x_n \}_{n\geq 1})=0$ or we note that
$ \{ x_n \}_{n\geq 1} \subset
\{\gamma F(\{ x_n\}_{n\geq 1})\,|\, \gamma\in [0,1] \}$, which gives
$$
\beta(\{x_n \}_{n\geq 1}) \leq \beta(\{\gamma F(\{ x_n\}_{n\geq 1})\,|\, \gamma\in [0,1] \})
= \beta(F(\{ x_n \}_{n\geq 1})) < \beta(\{ x_n \}_{n\geq 1}).
$$
But the latter case gives a contradiction and, therefore,
$(x_n)$ is relatively compact. Assume that $x_n\to x_0 \in \part U$. By continuity, it is clear that
 $x_0 = F(x_0)$, a contradiction with the admissibility of $I-F$,
 and the proof of (\ref{no-zeros-on-boundary-cond-2}) is finished.\\
\indent Thus one can put
\be\label{def-deg-condensing}
\deg (I-F,U):= \lim_{\lma\to 0^+} \deg(I-(1-\lma)F,U)
\ee
as $I-(1-\lma)F$ is a $k$-set contraction,
for any $\lma\in (0,1)$, and $(\ref{no-zeros-on-boundary-cond-2})$ holds. In view of the homotopy invariance property -- (D3), (\ref{def-deg-condensing}) is correct.
It can be easily verified that the extended degree admits
all the properties (D1)--(D4), meaning by {\em an admissible homotopy} a continuous map
$I-H:\o U\times [0,1]\to E$ such that
$0\not\in (I-H)(\part U\times [0,1])$ and, for any $\Omega\subset \o U$ with $\beta(\Omega)>0$,
$\beta(H(\Omega\times [0,1]))<\beta(\Omega).$

We shall now pass to the version of topological degree for perturbations of generators of $C_0$ semigroups of contractions. We say that a mapping $-A+F$ is {\em admissible} with respect to an open bounded $U\subset E$
if $A:D(A)\to E$ is a linear operator on a Banach space $E$ such that $-A$ generates a $C_0$
semigroup $S_A=\{ S_A(t):E \to E \}_{t\geq 0}$ of linear contractions, i.e. there exists $\omega>0$
such that
$$
\|S_A(t)\|\leq e^{-\omega t} \mbox{ for any } t\geq 0,
$$
$F:\o U \to E$ is a $k$-set contraction with a constant
$k\in [0,\omega)$ and $0\not\in (-A+F)(\part U\cap D(A))$.
By {\em an admissible homotopy} we shall understand a mapping
$(x,\lma)\mapsto -A(\lma)x + F(x,\lma)$ where the family of operators
$\{ A(\lma)\}_{\lma\in [0,1]}$ and continuous $F:\o U\times [0,1]\to E$ with open bounded $U\subset E$ satisfy
the following conditions
\begin{itemize}
\item the family of semigroups $\{S_{A(\lma)}\}_{\lma\in[0,1]} $ is uniformly contractive, i.e.
there is $\omega > 0$ such that $\|S_{A(\lma)} (t) \| \leq e^{-\omega t}$ for any $t\geq 0$
and $\lma\in [0,1]$,
\item the family $\{A(\lma)\}_{\lma\in [0,1]}$ is resolvent continuous,
that is if $\lma_n\to \lma_0$ (as $n\to\infty$) in $[0,1]$, then
$A (\lma_n) \stackrel{res}{\to} A (\lma_0)$ ({\em with respect to resolvent } or {\em graphs}), i.e.
for any $\nu > -\omega$ and $u\in E$, $R(\nu:-A(\lma_n))u \to R(\nu;-A(\lma_0))u$ as $n\to \infty$,
\item there exists $k\in [0,\omega)$ such that, for any $\Omega\subset \o U$,
$\beta (F(\Omega\times [0,1])) \leq k \beta (\Omega)$,
\item for any $\lma\in [0,1]$, $0\not\in (-A(\lma)+F(\cdot,\lma)) (\partial U\cap D(A(\lma)))$.
\end{itemize}
\begin{Rem}\label{after-admissible}{\em
The resolvent continuity condition for the family $\{ A(\lma)\}_{\lma\in [0,1]}$ of $C_0$ semigroups implies
that, for any $\lma_n \to \lma_0$ and $x\in E$, $S_{A(\lma_n)} (t) x \to S_{A(\lma_0)} (t)x$ uniformly
with respect to $t$ from bounded intervals (see \cite{Pazy}).  }
\end{Rem}
\begin{Lem} \label{resolvent-comp-property}
If $\{A(\lma) \}_{\lma\in [0,1]}$ and
$F:\o U\times [0,1]\to E$ determine an admissible homotopy,
then\\
(i) for any $t\geq 0$ and bounded $\Omega\subset E$,
$$
\beta \left(\bigcup_{\lma\in [0,1]} S_{A(\lma)}(t) \Omega \right) \leq
e^{-\omega t} \beta(\Omega);
$$
(ii) for $\nu \geq  -\omega$ and bounded $\Omega \subset E$,
$$
\beta \left(\bigcup_{\lma\in [0,1]} R(\nu:-A(\lma)) \Omega \right) \leq
\frac{1}{\omega+\nu}  \beta(\Omega)
$$
where $R(\nu:-A(\lma)):= (\nu I+ A(\lma))^{-1}$;\\
(iii) for any  $\nu \geq 0$,
the map $\Psi:\o U \times [0,1]\to E$, given by $\Psi(x,\lma):= R(\nu:-A(\lma))(\nu u +F(u,\lma))$,
is a set contraction with constant $\frac{k+\nu}{\omega + \nu}<1$ and $R(\nu:-A(\lma))(\nu x +F(x,\lma))\neq x$ for
any $(x,\lma) \in \part U\times [0,1]$.
\end{Lem}
{\bf Proof.} (i) Take an arbitrary $t>0$ and $\eps>0$. Let $\{ x_{k}\}_{k=1,\ldots,m_{\eps}}\subset E$ be
such that
$$
\Omega\subset \bigcup_{k=1}^{m_\eps} B \left( x_k, \beta(\Omega)+\eps/2\right).
$$
By the continuity property -- see Remark \ref{after-admissible},
the maps $[0,1]\ni \lma\mapsto S_{A(\lma)}(t)x_k\in E$, $k=1,\ldots,m_{\eps}$,
are uniformly continuous, therefore there
exist $\{ \lma_j \}_{j=1}^{n_\eps}$ and $\delta>0$ such that, for any
$k=1,\ldots,m_\eps$ and $j=1,\ldots,n_\eps$,
$$
\|S_{A(\lma)} (t)x_k - S_{A(\lma_j)}(t) x_k \| \leq \eps/2 \ \ \mbox { for any }
\lma\in B(\lma_j,\delta)
$$
and $\bigcup_{j=1}^{n_\eps} (\lma_j- \delta,\lma_j+\delta)\supset [0,1]$. Take any $x\in \Omega$ and
$\lma\in [0,1]$ and observe that, there are $k_x \in \{1,\ldots,m_\eps \}$ and
$j_{\lma}\in \{ 1,\ldots,n_\eps\}$, such that $\|x-x_{k_x}\| < \beta(\Omega)+\eps/2$ and
$|\lma - \lma_{j_\lma}|<\delta$. Further, by the contractivity,
$$
\begin{array}{ll}
\|S_{A(\lma)} (t)x - S_{A(\lma_{j_\lma})}(t)x_{k_x}\| & \leq
\|S_{A(\lma)} (t)x - S_{A(\lma)}(t)x_{k_x}\| + \\
& +\|S_{A(\lma)} (t)x_{k_x} - S_{A(\lma_{j_\lma})}(t)x_{k_x}\| \\
& \leq e^{-\omega t}\|x-x_{k_x} \| + \eps/2 \leq
r_\eps:= e^{-\omega t} \beta(\Omega) + \eps.
\end{array}
$$
Hence, the family $\{B(S_{A(\lma_j)}(t)x_k,r_\eps) \}_{j=1,\dots, n_\eps,\,
k=1,\ldots,m_\eps} $ covers the set $\bigcup\limits_{\lma\in [0,1]} S_{A(\lma)}(t)\Omega$,
which implies
$\beta\left(\bigcup\limits_{\lma\in [0,1]} S_{A(\lma)}(t)\Omega \right) \leq
r_\eps$. Since $\eps$ was arbitrary, the proof is finished.\\
\indent (ii) can be shown in the same manner as (i) with use of the inequality
$\|R(\nu:-A(\lma))\| \leq 1/(\omega+\nu)$ for any $\nu\geq -\omega$ and $\lma\in [0,1]$.\\
\indent (iii) Note that in view of (ii)
\begin{eqnarray*}
\beta  \left( \bigcup_{\lma\in [0,1]} R(\nu:-A(\lma))(\nu I+F)(\Omega\times \{ \lma \}) \right)
\mbox{ } \mbox{ } \mbox{ } \mbox{ } \mbox{ } \mbox{ } \mbox{ } \mbox{ } \mbox{ }
\mbox{ } \mbox{ } \mbox{ } \mbox{ } \mbox{ } \mbox{ } \mbox{ } \mbox{ } \mbox{ }
\mbox{ } \mbox{ } \mbox{ } \mbox{ } \mbox{ } \mbox{ } \mbox{ } \mbox{ } \mbox{ }
\mbox{ } \mbox{ } \mbox{ } \mbox{ } \mbox{ } \mbox{ } \mbox{ } \mbox{ } \mbox{ }
\mbox{ } \mbox{ } \mbox{ } \mbox{ } \mbox{ } \mbox{ } \mbox{ } \mbox{ } \mbox{ } \\
 \leq \beta \left( \bigcup_{\lma\in [0,1]} R(\nu:-A(\lma)) (\nu I+F)(\Omega\times [0,1])\right)
 \mbox{ } \mbox{ } \mbox{ } \mbox{ } \mbox{ } \mbox{ } \mbox{ } \mbox{ } \mbox{ }
\mbox{ } \mbox{ } \mbox{ } \mbox{ } \mbox{ } \mbox{ } \mbox{ } \mbox{ } \mbox{ }
\mbox{ } \mbox{ } \mbox{ } \mbox{ } \mbox{ } \mbox{ } \mbox{ } \mbox{ } \mbox{ } \\
 \leq\frac{1}{\omega+\nu} \beta ((\nu I+F)(\Omega\times [0,1])) \leq \frac{1}{\omega+\nu} (\nu\beta(\Omega)
+ k \beta(\Omega))= \left( \frac{k+\nu}{\omega+\nu} \right) \beta(\Omega).
\end{eqnarray*}
Finally, if $R(\nu:-A(\lma))(\nu I+F)(x,\lma)=x$ for some $x\in \part U$ and $\lma\in [0,1]$, then
$x\in D(A(\lma))$ and $-A(\lma)x+F(x,\lma)=0$, which contradicts the admissibility assumption.
\hfill $\square$

Therefore, if $-A+F$ is admissible with respect to an open bounded $U\subset E$, then, by Lemma
\ref{resolvent-comp-property}, we may put properly
\be\label{top-deg-A-F-def}
\Deg(-A+F, U):= \deg (I-R(\nu:-A)(\nu I+F),U)
\ee
where $\nu\geq 0$ is arbitrary. By the continuity of the resolvent and the homotopy invariance of $\deg$ it is clear that the definition (\ref{top-deg-A-F-def}) is independent of $\nu$.  It has also all the expected properties, i.e. the existence, additivity, homotopy
invariance and normalization. The latter property states that: {\em if $v_0 \in A(U\cap D(A))$,
then $\Deg(-A+v_0,U)=1.$} It is a straightforward conclusion from (\ref{top-deg-A-F-def}) and {\bf (D4)}.\\
\indent To define the topological degree of $-A+F$ in the situation $k=\omega$
it is sufficient to apply in (\ref{top-deg-A-F-def}) the topological degree for condensing vector fields.

\section{Degree formula for $k$-set contracting \\
and condensing fields}

Before we shall deal with the degree formula, we verify basic properties of
the translation along trajectories operator associated with the parameterized problem
$$
\begin{array}{ll}
\dot{u}(t) =-u(t) + G(u(t),\lma), \, t>0,\ \lma \in [0,1] \\
\end{array} \leqno{(P_{G,\lma})}
$$
where  $G: U \times [0,1]\to E$, with an open
bounded subset $U$ of a separable Banach space $E$, satisfy the following conditions\\
\noindent {\bf (G1)} \parbox[t]{120mm}{$G$ is locally Lipschitz in the first variable uniformly with respect to the second one, i.e. for any $x\in U$ there exist $\delta_x>0$ and $L_x>0$ such that, for any $x_1,x_2
\in B(x,\delta_x)\cap U$ and $\lma\in [0,1]$,
$$
\| G (x_1,\lma) - G (x_2, \lma)\|\leq L_x \|x_1-x_2\|;
$$}\\
\noindent {\bf (G2)} \parbox[t]{120mm}{$G(U\times [0,1])$ is bounded; }\\[1em]
\noindent {\bf (G3)} \parbox[t]{120mm}{$G$ is a $k$-set contraction, i.e.
there is $k\in [0,1)$ such that, for any bounded $\Omega\subset U$,
$$
\beta (G(\Omega\times [0,1])) \leq k \beta(\Omega).
$$ }

\noindent Let $V$ be any open  subset of $U$
such that $V \subset \o V + B(0,\rho)\subset U$ with some $\rho>0$.
For $t>0$, define $\Phi_t: \o V \to E$ by
$$
\Phi_t (x) := u(t;0, x,\lma)
$$
where $u(\cdot, 0,x,\lma)$ stands for the (local) solution of $(P_{G,\lma})$ with the initial condition
$u(0)=x$.
\begin{Prop} \label{prop-of-semiflow}
Under the above assumptions\\
{\em (i)} there exists $\o t>0$ such that, for any $x\in \o V$ and $\lma\in [0,1]$, there
exists a unique solution $u(\cdot;0,\o t, x,\lma)$ of $(P_{G,\lma})$ on $[0,\o t]$ starting at $x$
and the map
$$
\o V\times [0,1] \ni(x,\lma)\mapsto u(\cdot;0,\o t,x,\lma ) \in C([0,\o t],E);
$$
is continuous.\\
{\em (ii)} for any $t\in (0,\o t]$, $\Phi_t$ is a $k$-set
contraction, i.e. for any $\Omega\subset \o V$,
$$
\beta (\Phi_t (\Omega\times [0,1])) \leq e^{(k-1)t} \beta(\Omega).
$$
\end{Prop}
In the proof we shall use the following property of the measure of noncompactness.
\begin{Prop} \label{integral-prop-of-measure} {\em (see \cite{Deimling} or \cite{ObuKaZ})}
Let $E$ be a separable Banach space, $B\subset L^1([a,b],E)$ be countable and integrably bounded
(i.e. there exists $c\in L^1([a,b])$ such that $\|w(t)\|\leq c(t)$  for all $w\in B$ and a.e. $t\in [a,b]$)
and $\phi:[a,b]\to\R$ be given by $\phi (t):= \beta(\{u(t)|\, u\in B \}).$
Then $\phi\in L^1([a,b])$ and
$$
\beta \left(\left\{
\int_{a}^{b} u(\tau) \d \tau\,|\, u \in B\right\} \right)
\leq \int_{a}^{b} \phi(\tau) \d \tau.
$$
\end{Prop}
\noindent {\bf Proof of Proposition \ref{prop-of-semiflow}}.
(i) By the standard local existence and uniqueness theorem for ordinary differential equations
in a Banach space, we infer that, for any $x\in \o V$ and $\lma\in [0,1]$, there exists
a unique solution $u$ of $(P_{G,\lma})$ on some maximal interval $[0,\omega_x)$ with $\omega_x>0$.
By (G2), the exists $K>0$ such that $U \cup G(U\times [0,1]) \subset B(0,K)$. Hence
$$
\|u(t)-x\| \leq \int_{0}^{t} \|u(\tau)-G(u(\tau),\lma) \| \d \tau \leq 2K t \mbox{ for any }
t\in [0,\omega_x).
$$
This means that if $t\in [0, \rho /(4K)]$ and $t<\omega_x$, then $u(t)\in \o V + D(0,\rho/2) \subset U$, which
means that $\omega_x > \o t:= \rho/(4K)$. The continuity of $(x,\lma)\mapsto
u(\cdot;0,\o t,x,\lma)$ is a classical property of parameterized ordinary differential equations
in Banach spaces.
The proof of the part (i) is completed.\\
\indent (ii) Take any $\Omega\subset \o V$ and choose a countable set $\Omega_0 \subset \Omega$
such that $\Omega \subset\o {\Omega_0}$ and a countable and dense in $[0,1]$ subset
$\Lambda_0$. By the variation of constant formula, for any
$t\in [0,\o t]$, $x\in \o V$ and $\lma\in [0,1]$
$$
\Phi_t (x,\lma) = e^{-t} x + \int_{0}^{t} e^{-t+\tau} G(\Phi_\tau (x,\lma),\lma)\d\tau.
$$
This provides, in view of Proposition \ref{integral-prop-of-measure},
\begin{eqnarray*}
\beta (\Phi_t (\Omega_0\times \Lambda_0))  \leq e^{-t}\beta(\Omega_0) + \beta
\left(\left\{  \int_{0}^{t} e^{-t+\tau} G(\Phi_\tau (x,\lma),\lma)\d\tau\,|\,
 x\in \Omega_0,\,\lma\in \Lambda_0 \right\}\right)\\
   \leq e^{-t}\beta(\Omega_0) + \int_{0}^{t}
   \beta\left( e^{-t+\tau} G(\Phi_\tau(\Omega_0\times \Lambda_0)\times \Lambda_0) \right)\d\tau\\
\leq e^{-t}\beta(\Omega_0) + k \int_{0}^{t} e^{-t+\tau}  \beta (\Phi_\tau(\Omega_0\times \Lambda_0)) \d \tau.
\end{eqnarray*}
Hence, by the Gronwall inequality,
$e^{t}\beta(\Phi_t(\Omega_0\times \Lambda_0))\leq e^{kt} \beta(\Omega_0)$, which gives
\be\label{integral-prop-of-measure-1}
\beta(\Phi_t(\Omega_0\times \Lambda_0))\leq e^{(k-1)t} \beta(\Omega_0).
\ee
Finally, note that $\beta(\Omega)=\beta (\Omega_0)$, and,
by the continuity, $\Phi_t (\o{\Omega_0\times \Lambda_0}) \subset \o{\Phi_t (\Omega_0\times \Lambda_0)}$, which together with (\ref{integral-prop-of-measure-1}) gives
\begin{eqnarray*}
\beta(\Phi_t(\Omega\times [0,1])) = \beta (\Phi_t (\o{\Omega_0\times \Lambda_0}))
\leq  \beta(\o{\Phi_t (\Omega_0\times \Lambda_0)}) \\
= \beta (\Phi_t (\Omega_0\times \Lambda_0))
\leq e^{(k-1)t} \beta(\Omega_0) = e^{(k-1)t} \beta(\Omega).
\end{eqnarray*}
The proof of (ii) is completed. \hfill $\square$

We pass now to the first degree formula for the equation
$$
\dot u = - u + F(u) \leqno{(P_F)}
$$
where $F:U\to E$ is a locally Lipschitz bounded $k$-set contraction on
an open bounded set $U\subset E$.
\begin{Th} \label{Krasnosel-I}
Let  $F$ be as above and $V\subset U$ be an open set such that $\o V + B(0,\rho) \subset U$, for some $\rho>0$,
and $0\not\in (I-F)(\part V)$. Then there exists $\o t > 0$ such that, for any $t\in (0,\o t]$, \\
{\em (i)} the translation along trajectories operator for $(P_F)$,
$\Phi_t:\o V\to E$ is well-defined and $\Phi_t (x) \neq x$ for any $x\in \part V$;\\
{\em (ii)} $ \deg(I-F,V) = \deg(I-\Phi_t,V)$.
\end{Th}
{\bf Proof}. We reduce the proof to the case when $F$ is a compact map (see e.g.
\cite{Cwiszewski-1}).\\
\indent Due to Remark \ref{homotopy-compact}, there
is a locally Lipschitz compact map $F_L:\o V\to E$ homotopic to $F$ via
an admissible homotopy $H:\o V\times [0,1]\to E$, $H(x,\lma):=
(1-\lma)F(x)+\lma F_L(x)$. Then obviously
\be\label{Krasnosel-I-0}
\deg(I-F,V)=\deg(I-F_L,V).
\ee
Consider the parameterized problem given by
$$
\dot u = - u + H(u,\lma). \leqno{(P_{H,\lma})}
$$
Obviously, $H(\cdot,\lma)$ is locally Lipschitz and, by Proposition \ref{prop-of-semiflow} (i),
there exists $\o t_1>0$ such that, for any  $x\in \o V$ and $\lma\in [0,1]$, $(P_{H,\lma})$
admits a solution on $[0,\o t_1]$. Define $\Psi_t:\o V\times [0,1]\to E$ by
$$
\Psi_t (x,\lma):= u(t;0,\o t_1, x,\lma)
$$
where $u(\cdot;0,\o t_1, x,\lma)$ stands for the solution of $(P_{H,\lma})$ with the initial condition
$u(0)=x$.
By Proposition \ref{prop-of-semiflow} (ii), for any $\Omega\subset \o V$,
$$\beta(\Psi_t (\Omega\times [0,1])) \leq e^{(k-1)t} \beta(\Omega).$$

Now we claim that there exists $\o t_2\in (0,\o t_1]$ such that, for any $t\in (0,\o t_2]$,
\be\label{Krasnosel-I-1A}
\Psi_{t} (x,\lma)\neq x \mbox{ for any } (x,\lma)\in \part V\times [0,1].
\ee
Indeed, suppose to the contrary. Then there exist
sequences $t_{n}\rightarrow 0^+$, $(\lambda_{n})\subset [0,1]$ and
$(x_n )\subset\partial V$ such that $\Psi_{t_{n}}(x_n, \lma_n) = x_n$ for each $n\geq1$.
Note that this implies
$\Psi_{jt_n}(\lambda_n, x_n) = x_n$ for integer $j>0$ such that $j t_n \leq \o t_1$ and, in particular,
$\Psi_{\o t_1/2 + r_{n}}(x_n,\lma_n) = x_n,$ where
$r_n: = ([(\o t_1/2) / t_n ] + 1 ) t_n - \o t_1/2 < t_n \to 0^+$. Clearly, $0<r_n<t_n \to 0^+$ as
$n\to\infty$.
Since
$$\begin{array}{ll}
\|\Psi_{\o t_1/2+r_n}(x_n,\lma_n)-\Psi_{\o t_1/2} (x_n,\lma_n) \| & \leq
\displaystyle{\int\limits_{\o t_1/2}^{\o t_1/2+r_n} } \|-\Psi_{\tau}(x_n,\lma_n)+ F(\Psi_{\tau}(x_n,\lma_n))\| \d \tau\\
& \leq 2 K r_n  \to 0^+,
\end{array}$$
one gets
$$\begin{array}{ll}
\beta(\{x_{n}\}_{n\geq 1})    & =   \beta ( \{ \Psi_{\o t_1 / 2+r_{n}} ( x_n , \lma_n ) \}_{n\geq 1} )
                      =  \beta ( \{\Psi_{\o t_1 / 2} (x_n,\lma_n) \}_{n\geq 1}) \\
                       & \leq  e^{(k-1)\o t_1/2} \beta(\{x_{n}\}_{n\geq 1}),
\end{array}$$
which means that  $\beta(\{x_{n}\}_{n\geq 1}) = 0$. Hence, without loss of generality,
one may assume that $x_{n}\to x_{0}\in\part V$ and $\lma_n\to \lma_0\in [0,1]$.
Denote $u_n:= u(\cdot ; 0,\o t_1, x_n, \lambda_{n})$, for
$n\geq 0$, and observe that using the $t_n$-periodicity of $u_n$, for any $t\in [0,\o t_1]$, we get
\be \label{Krasnosel-I-1}
\| u_0 (t)-x_0 \|\leq \|u_0-u_n\|+ \|u_n(t)-u_n([t/t_{n}]t_{n})\|+\|x_n-x_{0}\|.
\ee
Furthermore, by Proposition \ref{prop-of-semiflow} (i), $u_n\to u_0$ in $C([0,\o t_1, E])$
and, this together with (\ref{Krasnosel-I-1}) provides $u_0 (t)= x_0$ for each $t\in [0,\o t_1]$,
i.e. $0 = - x_0 + H(x_0,\lma_0)$, a contradiction, since $I-H$ is an admissible
homotopy. Thus, the proof of (\ref{Krasnosel-I-1A}) is completed.
As an immediate conclusion of (\ref{Krasnosel-I-1A}), one obtains
\be\label{Krasnosel-I-2}
\deg(I-\Psi_t (\cdot,1),V) =\deg(I-\Psi_t (\cdot,0),V)\ \mbox{ for any } t\in (0,\o t_2].
\ee
\indent By the compact version of the theorem -- see \cite[Prop. 4.3]{Cwiszewski-1},
there exists $\o t \in (0,\o t_2]$
such that, for any $t\in (0,\o t]$,
\be\label{Krasnosel-I-3}
\deg ( I - F_L , V ) = \deg ( I-\Psi_{t} (\cdot, 1) , V ).
\ee
Finally, combining (\ref{Krasnosel-I-0}), (\ref{Krasnosel-I-3}) and (\ref{Krasnosel-I-2}),
we end the proof. \hfill   $\square$

Now let us extend Theorem \ref{Krasnosel-I} to the class of
condensing vector fields.
\begin{Th} \label{Krasnosel-II}
Let $F:U\to E$ be a locally Lipschitz bounded condensing (with
respect to the Hausdorff measure of noncompactness) map on an open
bounded set $U$ and $V\subset U$ be an open set such that $\o V +
B(0,\rho) \subset U$ and $0\not\in (I-F)(\part V)$. Then
there exists $\o t > 0$ such that, for any $t\in (0,\o t]$,\\
{\em (i)}  the translation along trajectories operator for $(P_F)$, $\Phi_t:\o V\to E$
is a well-defined condensing map and $\Phi_t (x) \neq x$ for any $x\in \part V$;\\
{\em (ii)} $ \deg(I-F,V) = \deg(I-\Phi_t,V)$.
\end{Th}
\begin{Rem}\label{after-Krasnosel-II} {\em
Under the assumptions of Theorem \ref{Krasnosel-II}, the thesis of Proposition
\ref{prop-of-semiflow} holds and the inequality in thesis (ii)
has the form $\beta(\Phi_t(\Omega\times [0,1]))\leq \beta(\Omega)$.}
\end{Rem}
{\bf Proof of Theorem \ref{Krasnosel-II}}.
Consider the parameterized family of equations
$$
\dot u = - u + (1-\lma)F(u) \ \ \ \lma\in [0,1]. \leqno{(P_{F,\lma})}
$$
By Remark \ref{after-Krasnosel-II}, there is $\o t_1>0$ such that
the map $\Phi_t:\o V\times [0,1]\to E$, $t\in [0,\o t_1]$, given by $\Phi_t(x,\lma):=u(t;0,\o t_1, x,\lma)$,
where $u(\cdot; 0,\o t_1, x,\lma)$ is the unique solution of $(P_{F,\lma})$ at time $0$ from $x$,
is well-defined, continuous and, for any $\Omega\subset \o V$,
$$
\beta(\Theta_t(\Omega\times [0,1]))\leq \beta(\Omega).
$$
We shall show that there exists $\o t\in [0, \o t_1]$, such that, for  $t\in [0,\o t]$, $I-\Theta_t$
is an admissible homotopy in the sense of the topological degree for condensing maps.
First, we need to prove
that for any $t\in (0,\o t_1]$ and $\Omega\subset \o V$ with $\beta(\Omega)>0$,
\be\label{Krasnosel-II-1}
\beta(\Theta_t(\Omega\times [0,1]))< \beta(\Omega).
\ee
Take any fixed $t\in (0,\o t_1]$, $\Omega\subset  \o V$ with $\beta(\Omega)>0$, $\Omega_0\subset \Omega$ a countable set such that $\Omega \subset \o{\Omega_0}$ and $\Lambda_0 \subset [0,1]$  a countable dense subset of $[0,1]$.
Then, obviously, for any $t\in [0, t]$, by Proposition \ref{integral-prop-of-measure},
\begin{eqnarray*} \label{Krasnosel-II-2-bis}
\nonumber \beta(\Theta_t(\Omega_0 \times \Lambda_0 )) \leq e^{-t}\beta(\Omega_0)
+ \int_{0}^{t} e^{-t+\tau}
\beta\left(\bigcup_{\gamma\in [0,1]\cap\Q} \gamma F(\Theta_\tau(\Omega_0\times \Lambda_0 ))\right) \d\tau\\
\leq e^{-t}\beta(\Omega_0)\nonumber
+ \int_{0}^{t} e^{-t+\tau}
\beta\left(F(\Theta_\tau(\Omega_0\times \Lambda_0 ))\right) \d\tau. \nonumber
\end{eqnarray*}
Clearly, either $\beta(\Theta_\tau(\Omega_0\times \Lambda_0))=0$, for a.e. $\tau\in [0,t]$
and (\ref{Krasnosel-II-1}) follows clearly or there exists $J\subset [0,t]$ of positive Lebesgue measure such that $\beta(\Theta_\tau(\Omega_0\times \Lambda_0))>0$ for any $\tau\in J$. Hence
\begin{eqnarray*}
\beta(\Theta_t(\Omega_0 \times \Lambda_0 )) <  e^{-t}\beta(\Omega_0)+ \int_{0}^{t} e^{-t+\tau}
\beta\left(\Theta_\tau(\Omega_0\times \Lambda_0 )\right) \d\tau\\
\leq e^{-t}\beta(\Omega_0)+ \int_{0}^{t} e^{-t+\tau}
\beta(\Omega_0) \d\tau = \beta(\Omega_0),
\end{eqnarray*}
which implies (\ref{Krasnosel-II-1}).\\
\indent Next, one has to show that there exist
$\lma_1 \in (0,1]$ and $\o t_2\in (0,\o t_1]$ such that, for any $t\in (0,\o t_2]$,
$$
\Theta_t (x,\lma)\neq x \ \ \mbox{ for } (x,\lma)\in \part V \times [0,\lma_1].
$$
If we suppose to the contrary that, then there exist
$(x_n)\subset \part V$, $\lma_n\to 0^+$ and $t_n\to 0^+$ such that
$\Theta_{t_n} (x_n,\lma_n)=x_n$ for each $n\geq 1$. This clearly implies
that $\Theta_{jt_n}(x_n, \lma_n)=x_n$ for any integer $j>0$ such that $j t_n \leq \o t_1$,
which gives, for large $n$, $\Theta_{\o t_1/2 + r_n} (x_n, \lma_n) = x_n$ where
$r_n:=([(\o t_1/2)/t_n]+1) t_n - (\o t_1/2)$. Clearly $0<r_n<t_n \to 0^+$ and
\begin{eqnarray*}
\|\Theta_{\o t_1/2+r_n} (x_n,\lma_n)-\Theta_{\o t_1/2}(x_n,\lma_n)\|
\leq \int\limits_{\o t_1/2}^{\o t_1/2 + r_n} \|-\Theta_\tau(x_n,\lma_n)+
(1-\lma_n)F(\Theta_\tau(x_n,\lma_n))\| \d \tau\\
\leq 2 K r_n\to 0^+
\end{eqnarray*}
where $K>0$ is such that $U\cup F(U) \subset B(0,K)$.
Therefore either $\beta(\{ x_n \}_{n\geq 1})=0$ or
\begin{eqnarray*}
\beta(\{ x_n\}_{n\geq 1}) = \beta(\{ \Theta_{\o t_1/2+r_n}(x_n,\lma_n) \}_{n\geq 1})
= \beta(\{\Theta_{\o t_1/2}(x_n,\lma_n) \}_{n\geq 1}) < \beta(\{x_n \}_{n\geq 1}).
\end{eqnarray*}
The latter inequality gives a contradiction and we infer that
$(x_n)$ is relatively compact. Without loss of generality we may assume that
$x_n\to x_0\in \part V$. If one puts $u_n:= u(\cdot;0,\o t_1,x_n,\lma_n)$, then,
by use of Proposition \ref{prop-of-semiflow} (i),
$u_n\to u_0$ in $C([0,\o t_1],E)$ and, by the same argument as in  (\ref{Krasnosel-I-1}),
$u_0 \equiv x_0$, that is  $0= - x_0+F(x_0)$, a contradiction.\\
\indent Finally, according to (\ref{def-deg-condensing}), there is $\lma_0 \in (0,\lma_1]$ such that $\deg(I-F,V) = \deg(I-(1-\lma_0)F,V)$ and, in view of Theorem \ref{Krasnosel-I}, there exists $\o t \in (0,\o t_2]$
such that, for any $t\in (0, \o t]$,
$\deg (I-F, V) = \deg ( I-(1-\lma_0)F , V ) = \deg( I-\Theta_t ( \cdot , \lma_0) , V ).$
This together with  the homotopy invariance for the topological degree of condensing
vector fields yields $\deg (I-F,V ) = \deg(I-\Theta_t (\cdot,0), V)=\deg(I-\Phi_t,V).$ \hfill $\square$

\section{Differential equations governed by \\
perturbations
of generators of  contractive \\
$C_0$ semigroups} \label{section-3}

\noindent We provide some general properties of solution operators
generated by perturbations of generators of $C_0$ semigroups.
They enable us to apply homotopy arguments in the next section.

Suppose that $\{ A(\lma):D(A(\lma))\to E\}_{\lma\in [0,1]}$ is a family of densely defined linear operators such that

\noindent {\bf (A1)} \parbox[t]{140mm}{ for any $\lma\in [0,1]$,  $-A(\lma)$ is an infinitesimal generator of a $C_0$ semigroup
$\{ S_{A(\lma)}(t) \}_{t\geq 0}$ of bounded linear operators on $E$ and there exists $\omega>0$ such that
$\|S_{A(\lma)} (t)\| \leq e^{-\omega t}$ for any $t\geq 0$ and $\lma\in [0,1]$;}\\[1em]
\noindent {\bf (A2)} \parbox[t]{140mm}{ if $\lma_n\to \lma_0$ (as $n\to\infty$) in $[0,1]$, then
$A (\lma_n) \stackrel{res}{\to} A (\lma_0)$ ({\em with respect to graphs } or {\em resolvent}), i.e.
for any $\nu > -\omega$ and $u\in E$, $R(\nu:-A(\lma_n))u \to R(\nu;-A(\lma_0))u$,}\\[1em]
\noindent and that $F:[0,T]\times E\times [0,1] \to E$, where $T>0$,
is a continuous mapping satisfying the following conditions

\noindent {\bf (F1)} \parbox[t]{140mm}{$F$ is locally Lipschitz in the second variable
uniformly with respect to the others, i.e. for any $x \in E$ there exist $\delta_x>0$ and $L_x>0$ such that,
for any $x_1,x_2 \in B(x,\delta_x)$, $t\in [0,T]$ and  $\lma\in [0,1]$,
$\|F(t,x_1,\lma)-F(t,x_2, \lma)\|\leq L_x \|x_1-x_2\|;$}\\
\noindent {\bf (F2)} \parbox[t]{140mm}{$F$ has a (uniform) linear growth, i.e.
there exists $c>0$ such that $\|F(t,x,\lma)\|\leq c(1+\|x\|)$ for any $x\in E$,
$t\in [0,T]$ and $\lma\in [0,1]$; }\\[1em]
\noindent {\bf (F3)} \parbox[t]{140mm}{there is $k\geq 0$ such that, for any bounded $\Omega\subset E$,
$\beta(F([0,T]\times \Omega\times [0,1])) \leq k\beta(\Omega).$}\\

\noindent Below we collect some basic facts concerning mild solutions (see \cite{Pazy}) for the family of differential equations
$$
\dot u(t) = -A(\lma) u(t) + F(t,u(t),\lma) \ \mbox { on } [0,T]. \leqno{(P_{A,F,\lma})}
$$
\begin{Prop} \label{semiflow-props-A-F}$\mbox{ } $\\
{\em (i)} \parbox[t]{140mm}{{\em (Existence)} For any $x\in E$ and $\lma \in [0,1]$,
the equation $(P_{A,F,\lma})$ with the initial condition $u(0)=x$ has a unique mild solution
$u(\cdot;0,T,x,\lma)$. Moreover, if $T=+\infty$, then the assertion holds, too.}\\[1em]
{\em (ii)} \parbox[t]{140mm}{{\em (Continuity)} If $( x_n, \lma_n )\to ( x_0 , \lma_0 )$
in $E\times [0,1]$, then
$$
u(\cdot;0,T,x_n,\lma_n) \to u(\cdot;0,T,x_0,\lma_0) \ \ \ \mbox{ in } C([0,T],E).
$$}\\[1em]
{\em (iii)} \parbox[t]{140mm}{{\em (Compactness)} If $\Phi_t:E\times [0,1] \to E$,
for $t\in [0,T]$, is given by
$$  
\Phi_t (x,\lma):= u(t;0,T,x,\lma),$$
then, for any bounded $\Omega\subset E$,
$$\beta \left( \Phi_t (\Omega\times [0,1]) \right) \leq  e^{(k-\omega) t} \beta(\Omega).$$}
\end{Prop}
\begin{Rem} \label{semiflow-gronwall}
{\em (a) If $u:[0,T]\to E$ is a mild solution of $(P_{A,F,\lma})$, for some fixed $\lma\in[0,1]$,
with the initial condition $u(0)=x$ for some $x\in E$, then
$\|u(t)\| \leq \|S_{A(\lma)} (t) x \|+\int_{0}^{t} \|S_{A(\lma)}(t-\tau) F(\tau,u(\tau), \lma) \| \d\tau
 \leq e^{-\omega t} \|x\| + \int_{0}^{t} e^{-\omega(t-\tau)}c (1+\|u(\tau)\|) \d \tau$, for any $t\in [0,T]$,
and, in view of the Gronwall inequality, $\|u(t)\|\leq (\|x\|+ cTe^{\omega T})e^{cT}$,
for $t\in [0,T]$.\\
\indent (b) It follows from (a) that if $\Omega \subset E$ is bounded,
then $\{ \Phi_t (x,\lma) \, | \, (x,\lma) \in \Omega \times [0,1], \, t\in [0,T] \}$
is bounded.}
\end{Rem}
{\bf Proof}. The proof of (i) follows from  the arguments of \cite[Ch. 6, Th. 2.2]{Pazy},
the local existence theorem and Remark \ref{semiflow-gronwall}.\\
\indent To justify (ii) and (iii), observe that, by the definition of mild solution,
$$
\Phi_t (x,\lma) = S_{A(\lma)} (t)x + \int_{0}^{t} S_{A(\lma)}(t-\tau)
F(\tau,\Phi_\tau (x,\lma), \lma) \d\tau,
$$
for any $x\in \Omega$, $\lma\in [0,1]$ and $t\in [0,T]$. In view of Remark \ref{semiflow-gronwall},
the set $\Phi_t (\Omega\times [0,1])$ is bounded. Hence, acting with the measure of noncompactness, one obtains
\begin{eqnarray}\label{semiflow-props-A-F-0}
\ \ \ \ \beta\left(\Phi_t (\Omega\times [0,1])\right)
\leq
\beta\left(\bigcup\limits_{\lma\in [0,1]} S_{A(\lma)}(t)\Omega \right)+  \\
+\beta\left(\left\{ \int\limits_{0}^{t} w_{x,\lma,t}(\tau)\d\tau| x\in\Omega,
\lma\in[0,1]\right\}\right) \nonumber
\end{eqnarray}
where $w_{x,\lma,t}(\tau):=S_{A(\lma)}(t-\tau)F(\tau,\Phi_\tau(x,\lma),\lma)$.\\
\indent Suppose $\Omega_0$ is a countable subset of $\Omega$ such that $\Omega \subset \o{\Omega_0}$
and $\Lambda_0$ is a countable dense subset of $[0,1]$.
Then, in view of Proposition \ref{integral-prop-of-measure}, for any $t\in [0,T]$,
\begin{eqnarray*}
\beta\left( \left\{ \int_{0}^{t} w_{x,\lma,t}(\tau)\d\tau\,|\, x\in\Omega_0,\,
\lma\in\Lambda_0\right\}\right) \leq
\int_{0}^{t} \beta \left(\{w_{x,\lma,t}(\tau)
\,|\, x\in\Omega_0,\, \lma\in \Lambda_0  \}\right) \d\tau\\
\leq \int_{0}^{t} \beta \left( \bigcup_{\lma\in\Lambda_0}
S_{A(\lma)} (t-\tau) [F(\{\tau\}\times \Phi_\tau(\Omega_0\times\Lambda_0)\times\Lambda_0)] \right) \d\tau.
\end{eqnarray*}
and, in view of Lemma \ref{resolvent-comp-property} and (F3), one has
\begin{eqnarray}\label{semiflow-props-A-F-2} \nonumber
\beta\left( \left\{ \int_{0}^{t} w_{x,\lma,t} (\tau)\d\tau\,|\, x\in\Omega_0,\,
\lma\in\Lambda_0\right\} \right)   \mbox{ } \mbox{ } \mbox{ } \mbox{ } \mbox{ } \mbox{ } \mbox{ } \mbox{ } \mbox{ } \mbox{ } \mbox{ } \mbox{ } \mbox{ } \mbox{ } \mbox{ } \mbox{ } \mbox{ } \mbox{ } \mbox{ } \mbox{ } \\
\mbox{ } \mbox{ } \mbox{ } \mbox{ } \mbox{ } \mbox{ } \mbox{ } \mbox{ } \leq \int_{0}^{t}  e^{-\omega(t-\tau)} \beta(F([0,T]\times \Phi_\tau(\Omega_0\times\Lambda_0)\times\Lambda_0))\d\tau\\ \nonumber
\leq \int_{0}^{t}  k e^{-\omega(t-\tau)} \beta (\Phi_\tau(\Omega_0\times\Lambda_0))\d\tau.
\end{eqnarray}
Hence, by (\ref{semiflow-props-A-F-0}) and (\ref{semiflow-props-A-F-2}), one has
$$
e^{\omega t}\beta(\Phi_t (\Omega_0\times \Lambda_0)) \leq
\beta(\Omega_0) + \int_{0}^{t} k e^{\omega \tau}
\beta (\Phi_\tau(\Omega_0\times\Lambda_0))\d\tau,
$$
which, by the Gronwall inequality, yields
\be\label{semiflow-props-A-F-3}
\beta(\Phi_t(\Omega_0 \times \Lambda_0)) \leq e^{(k-\omega)t} \beta(\Omega_0).
\ee
Before we conclude (iii) for $\Omega$, we shall prove (ii). To this end take any sequence $(x_n, \lma_n)\subset E\times[0,1]$ with $(x_n,\lma_n)\to (x_0,\lma_0)$. In view of (\ref{semiflow-props-A-F-3}),
for any $t\in [0,T]$,
$$
\beta(\{ u (t;0,T,x_n, \lma_n)\}_{n\geq 1}) \leq  e^{(k-\omega)t} \beta(\{ x_n\}_{n\geq 1}) = 0,
$$
which together with (F2) and \cite[Prop. 2.7]{Cwiszewski-1}, implies the relative compactness of
$\{ u(\cdot;0,T, x_n, \lma_n)\}_{n\geq 1}$ in $C([0,T],E)$. Hence, any subsequence of
$(u(\cdot;0,T, x_n, \lma_n) )$, contains a convergent subsequence. Without loss of generality we may assume that the sequence of solutions $u_n:=u(\cdot; 0,T, x_n, \lma_n)$ converges in $C([0,T],E)$ to some $u_0$.
Since, for any $n\geq 1$, $u_n(t)=S_{A(\lma_n)} (t)x_n + \int_{0}^{t}
S_{A(\lma_n)} (t-\tau) F(\tau,u_n(\tau),\lma_n) \d \tau$, passing to the limit and using Remark
\ref{after-admissible}, one obtains
$u_0(t)=S_{A(\lma_0)}(t)x_0 + \int_{0}^{t} S_{A(\lma_0)} (t-\tau) F(\tau,u_0(\tau),\lma_0) \d \tau$,
i.e. $u_0 = u(\cdot;0,T,x_0,\lma_0)$. Thus, since we have shown that any subsequence of
$(u(\cdot;0,T,x_n,\lma_n))$ contains a subsequence converging to $u(\cdot;0,T,x_0,\lma_0)$, the assertion of (ii) is proved.\\
\indent Return now to (iii) for $\Omega$. Clearly, $\beta(\Omega)=\beta(\Omega_0)$ and, by the continuity from (ii), for any $t\in [0,T]$, $\Phi_t (\o{\Omega_0\times \Lambda_0}) \subset \o{\Phi_t (\Omega_0\times \Lambda_0)}$,
therefore $\beta(\Phi_t(\Omega\times [0,1])) \leq \beta(\o{\Phi_t (\Omega_0\times \Lambda_0)})
= \beta(\Phi_t (\Omega_0\times \Lambda_0)) \leq e^{(k-\omega)t} \beta(\Omega_0)=
e^{ (k-\omega) t } \beta(\Omega)$. \hfill $\square$

\begin{Rem} {\em
In the above theorem, it is possible to relax or get rid of the separability assumption on $E$
at cost of stronger assumptions on $k$ or geometry of $E$ -- see e.g. \cite{BaKr2}.
}\end{Rem}

\section{Krasnosel'skii type degree formula}
Now we are concerned with the differential equation
$$
\dot{u}(t) =-Au(t)+F(u(t)), \,\,t>0 \leqno{(P_{A,F})}
$$
where $A:D(A)\to E$ is a densely defined linear operator on a separable Banach space such that
$-A$ is an infinitesimal generator of a $C_0$ semigroup $\{S_A(t)\}_{t\geq 0}$
of bounded linear operators and $F:E\rightarrow E$ is a locally Lipschitz map
and suppose the following properties hold

\noindent {\bf (H1)} \parbox[t]{140mm}{ there exists $\omega>0$ such that
$\|S_A(t)\|_{{\cal L}(E,E)} \leq e^{-\omega t}$ for any $t\geq 0$;}\\ [1em]
\noindent {\bf (H2)} \parbox[t]{140mm}{ $F$ has {\em sublinear growth}, i.e. there exists $c>0$
such that $\|F(x)\|\leq c(1+\|x\|)$ for any $x\in E$;}\\ [1em]
\noindent {\bf (H3)} \parbox[t]{140mm}{ there is
$k\in [0,\omega)$ such that, for any bounded $\Omega\subset E$,
$\beta(F(\Omega)) \leq k\beta(\Omega).$}

By $\Phi_t:E\to E$, for $t>0$, denote the translation along trajectories operator given by
$$
\Phi_t (x):= u (t;0,x)
$$
where $u(\cdot;0,x)$ is a {\em mild solution } of $(P_{ A , F })$ with the initial value
condition $u(0)=x$. It is well-defined due to Proposition \ref{semiflow-props-A-F}.

Let us now state one of the main results of the paper, i.e. the Krasnosel'skii type formula
for $(P_{A,F})$.
\begin{Th} \label{degree-formula-II}
Let $A$ and $F$ be as above. Then, if $V\subset E$ is open bounded and $0 \not\in (-A +
F)(\partial V\cap D(A))$, then there is $\o t >0$ such  that for
any $t\in (0, \o t]$ and $x\in \part V$,  $\Phi_t (x)\neq x$,
and
$$
\deg(I-\Phi_t,V)= \Deg(-A+F,V).
$$
\end{Th}
\begin{Lem} \label{degree-formula-II-lemma}
Assume that $\{ A(\lma) \}_{\lma\in [0,1]}$ and $F:E\times [0,1]\to E$ satisfy
assumptions $(A1)$--$(A2)$ and $(F1)$--$(F3)$ from Section \ref{section-3} and $k\in [0,\omega)$.
Let $\Phi_t: E \times [0,1] \to E$, $t\geq 0$, be given by $\Phi_t (x,\lma):= u(t;0,x,\lma)$
where $u(\cdot;0,x,\lma)$ is a unique mild solution of
$$
\left\{ \begin{array}{l}
\dot u = -A(\lma)u + F(u,\lma) \ \mbox{ on } [0,+\infty) \\
u (0) = x.
\end{array} \right.
$$
If $V\subset E$ is an open bounded set such that $0\not\in (-A+F)(\part V\cap D(A)\times [0,1])$,
then there exists $\o t>0$ such that, for any $t\in (0,\o t]$,
$$
\Phi_t (x, \lma) \neq x \ \mbox{ for any } (x,\lma)\in \part V\times [0,1].
$$
\end{Lem}
\noindent {\bf Proof}.
Suppose to the contrary that there exist $t_n\to 0^+$, $(x_n)\subset \part V$
and $(\lma_n)\subset [0,1]$ such that
$\Phi_{t_n} (x_n,\lma_n) = x_n$ for any $n\geq 1$. Then, for any $N\geq 1$,
\be\label{degree-formula-II-2}
\{ x_n \}_{n\geq 1} \subset \{ \Phi_{N+s_n(N)}(x_n,\lma_n) \}_{n\geq 1}
=  \{\Phi_N (\Phi_{s_n(N)} (x_n,\lma_n),\lma_n) \}_{n\geq 1}
\ee
where $s_n (N) := ([N / t_n]+1) t_n-N$. Clearly $0<s_n(N)\leq t_n$.
In view of Remark \ref{semiflow-gronwall}, there exists $R>0$ such that, for any
$N\geq 1$, $\{ \Phi_{s_n (N)} (x_n,\lma_n) \}_{n\geq 1} \subset B(0,R)$.
Hence, by use of (\ref{degree-formula-II-2}) and Proposition \ref{semiflow-props-A-F} (iii),
one has
\begin{eqnarray*}
\beta(\{ x_n \}_{n\geq 1})
\leq \beta ( \Phi_N ( \{ \Phi_{s_n (N)} ( x_n , \lma_n ) \}_{n\geq 1} \times [0,1]) ) \\
\leq e^{(k-\omega) N} \beta(\{ \Phi_{s_n (N)} (x_n,\lma_n) \}_{n\geq 1} )
\leq R e^{(k-\omega) N}.
\end{eqnarray*}
Passing to the limit with $N\to\infty$, one gets $\beta(\{ x_n\}_{n\geq 1}) = 0$. Hence,
without loss of generality, one may assume that $x_n\to x_0 \in \part V$ and
$\lma_n\to \lma_0$ in $[0,1]$.
Let $u_n:[0,1]\to E$ be given by $u_n (t):= u(t;0,x_n,\lma_n)$, for any $t\in [0,1]$ and $n\geq 0$.
Then, by Proposition \ref{semiflow-props-A-F} (ii), $u_n\to u_0$ in $C([0,1],E)$,
and, consequently, using the $t_n$-periodicity of $u_n$ and the equicontinuity of $\{ u_n \}_{n\geq 1}$,
one obtains, for any $t\in [0,1]$,
$$
\|u_0(t)-x_0\| \leq \|u_0-u_n \| + \|u_n(t) - u_n ([t/t_n]t_n)\|
+ \| x_n - x_0 \| \to 0 \ \ \mbox{ as } n\to\infty,
$$
i.e. $u_0 \equiv x_0$ on $[0,1]$.
By the definition of mild solution, one has, for any $n\geq 1$ and $t\in [0,1]$,
$u_n(t) = S_{A(\lma_n)}(t) x_n + \int_{0}^{t} S_{A(\lma_n)}(t-\tau)F(u_n(\tau),\lma_n) \d \tau.$
Passing to the limit, for any $t\in [0,1]$,
$$
x_0=S_{A(\lma_0)}(t) x_0 + \int_{0}^{t} S_{ A(\lma_0)}(t-\tau) F(x_0,\lma_0)\d\tau,
$$
which yields
$$
\lim_{t\to 0^+} \frac{1}{t}(x_0 - S_{ A(\lma_0)}(t)x_0) =
\lim_{t\to 0^+} \frac{1}{t} \int_{0}^{t} S_{A(\lma_0)}(t-\tau) F(x_0,\lma_0)\d\tau = F(x_0,\lma_0),
$$
i.e. $x_0\in D(A (\lma_0))$ and $0=-A(\lma_0)x_0 + F(x_0,\lma_0)$, a contradiction. \hfill$\square$

\noindent {\bf Proof of Theorem \ref{degree-formula-II}}.\\
\noindent {\bf Step 1}.
Assume first that $F:E\to E$ is completely continuous and define
the family $\{\widetilde A(\lma) \}_{\lma\in [0,1]}$
and the map $\widetilde F:E\times [0,1]\to E$ by
$$
\widetilde A(\lma):=(1-\lma)I+\lma A \,\,\, \mbox { for } \lma \in [0,1],
$$
$$
\widetilde F(x,\lma):=[\lma I  + (1-\lma)A^{-1}] F(x) \,\,\, \mbox{ for }
x\in E,\, \lma\in [0,1].
$$
By a direct verification (with use of the resolvent identity and limit property) it can be seen that
$\widetilde A (\lma_n) \stackrel{res}{\to} \widetilde A(\lma_0)$ whenever
$\lma_n\to \lma_0$ (as $n\to \infty$) and that
$\| S_{\widetilde A(\lma)}(t) \|_{{\cal L}(E,E)} =  \| e^{-(1-\lma)t} S_{\lma A}(t) \|_{{\cal L}(E,E)} =
\| e^{-(1-\lma)t}S_{A}(\lma t) \|_{{\cal L}(E,E)} \leq    e^{-(1-\lma)t - \lma\omega t}  \leq
   e^{-\widetilde \omega t}$ where $\widetilde \omega := \min \{1 , \omega \}$, i.e.
$\{ \widetilde A(\lma)\}_{\lma\in [0,1]}$ satisfies conditions $(A1)$--$(A2)$.
It is obvious that $\widetilde F$ satisfies $(F1)$--$(F2)$ and,
 for any bounded $\Omega \subset E$,
$\widetilde F (\Omega\times [0,1]) \subset \conv [F(\Omega )\cup A^{-1} F(\Omega )]$
is relatively compact. Moreover,
\be\label{degree-formula-II-A}
-\widetilde A(\lma)x + \widetilde F (x,\lma) \neq 0 \ \ \ \mbox{ for any } (x,\lma)\in \part V\cap D(A(\lma))\times [0,1].
\ee
Indeed for $\lma=0$ it is obvious and if there exists $(x_0,\lma_0)\in \part V\cap D(A)\times (0,1]$ such that
$-\widetilde A(\lma_0)x_0+\widetilde F(x_0,\lma_0)=0$, then
\begin{eqnarray*}
[(1-\lma_0)I+\lma_0 A]x_0 = [\lma_0 I  + (1-\lma_0) A^{-1}] F (x_0)\\
x_0 = ((\lma_{0}^{-1}(1-\lma_0))I + A)^{-1}[I +\lma_{0}^{-1}(1-\lma_0) R(0:-A)] F (x_0)\\
x_0 = R(\mu_0;-A)(I + \mu_0 R(0;-A))F(x_0)
\end{eqnarray*}
where $\mu_0:=\lma_{0}^{-1}(1-\lma_0)$. By the resolvent identity, $x_0 = R(0;-A) F(x_0)$, i.e.
$-Ax_0+F(x_0)=0$, a contradiction with the assumption, which completes the proof of (\ref{degree-formula-II-A}).\\
\indent In view of Proposition \ref{semiflow-props-A-F}, for any $t\in (0,1]$,
the map $\Psi_t:E\times [0,1]\to E$, given by
$$\Psi_t(x,\lma):=u(t;0,1,x,\lma)$$
where $u(\cdot;0, 1, x,\lma)$ is a unique mild solution of
$$\left\{\begin{array}{l}
\dot u (t) = -\widetilde A(\lma) u(t) + \widetilde F(u(t), \lma) \ \mbox{ on } [0,1]\\
u(0)=x,
\end{array}\right.$$
is continuous and, for any bounded $\Omega\subset E$,
\be\label{degree-formula-II-1}
\beta(\Psi_t ( \Omega \times [0,1] )) \leq e^{ -\widetilde \omega t} \beta(\Omega).
\ee
By (\ref{degree-formula-II-A}) and Lemma \ref{degree-formula-II-lemma}, there exists $\o t_1>0$
such that, for any $t\in (0,\o t_1]$, $\Psi_t (x,\lma)\neq x$, for any
$(x,\lma)\in \part V\times [0,1]$. Hence $I-\Psi_t$ is an admissible homotopy in the degree
theory and, by use of the homotopy invariance, for $t\in (0,\o t_1]$,
\be\label{degree-formula-II-4}
\deg (I-\Phi_t,V) = \deg (I-\Psi_t(\cdot,1),V) = \deg (I-\Psi_t(\cdot,0),V).\ee
In view of Theorem \ref{Krasnosel-I}, there exists $\o t_2\in (0,\o t_1]$, such that,
for any $t\in (0,\o t_2]$,
$$\deg(I-\Psi_t (\cdot,0),V) = \deg(I-R(0;-A)F,V),$$
which together with (\ref{degree-formula-II-4}) and the definition of the degree
completes the proof for completely continuous $F$.

\noindent {\bf Step 2}. Now suppose that $F$ satisfies (H3)
with $k\in [0,\omega)$. Obviously, $A^{-1}F:E \to E$ is a $k$-set contraction with the
constant $k/\omega < 1$. Let $C\subset E$ be a fundamental set for $A^{-1}F|_{\o V}$
(see Remark \ref{construct-fundset}) and $F_0:E\to E$ be a continuous extension of
$F |_{C\cap \o V}$ such that $F_0(E)\subset \o\conv F(C\cap\o V)$ -- existing by the Dugundji
extension theorem (see e.g. \cite{Dugundji-Granas}).
Furthermore let $F_L:E\to E$ be a locally Lipschitz map such that, for any $x\in E$,
\be \label{degree-formula-III-step-4-1}
\|F_L(x)-F_0(x)\|\leq \rho\omega/2,
\ee
with $\rho:= \inf\{ \|x-A^{-1}F(x)\|\,|\, x\in \part V \}>0$, and
$F_L(E)\subset \o\conv F_0(E)$ -- existing by the Lasota-Yorke approximation theorem.
Note that since $F_L (E) \subset \o\conv F_0(E)\subset \o\conv F(C\cap \o V)$, $F_L$ is compact.
Observe also that
$A^{-1}(F_L(E)) \subset A^{-1}(\o\conv F(C\cap \o V)) \subset
\o{A^{-1} (\conv F(C\cap \o V)) }  = \o\conv A^{-1}F(C\cap \o V) \subset
\o\conv C = C$.\\
\indent Now define $\o F: E\times [0,1]\to E$ by $\o F(x,\lma):=(1-\lma)F(x)+\lma F_L(x)$.
It is clear that, for any bounded $\Omega\subset E$, $\beta(\o F(\Omega\times [0,1]))
\leq k \beta(\Omega)$. Moreover, note that
\be\label{degree-formula-III-step-4-2}
-Ax+ \o F(x,\lma)\neq 0 \mbox{ for any } x\in \part V\cap D(A) \mbox{ and  } \lma\in [0,1].
\ee
Indeed, if $-Ax+\o F(x,\lma)=0$ for some $(x,\lma)\in \part V\cap D(A)\times [0,1]$, then
$x = A^{-1} \o F(x,\lma) \in A^{-1} (\conv [F(x)\cup F_L(\o V)])
 \subset \conv [A^{-1} F(x)\cup A^{-1} F_L(\o V)]\subset \conv[A^{-1}F(x) \cup C]$,
 and, since $C$ is fundamental for $A^{-1}F|_{\o V}$, we infer that $x\in C$; hence,
 by use of (\ref{degree-formula-III-step-4-1}),
\begin{eqnarray*}
\|x-A^{-1} F(x)\| =\lma \|A^{-1}F_L(x)-A^{-1}F(x)\| \leq (1/\omega) \|F_L(x)-F_0 (x)\| \\
                  \leq (1/\omega) \rho\omega/2 \leq \rho/2,
\end{eqnarray*}
which contradicts the definition of $\rho$ and proves (\ref{degree-formula-III-step-4-2}).\\
\indent Using the homotopy invariance of the degree for $-A+\o F$ we obtain
\be\label{degree-formula-III-step-4-3}
\Deg(-A+F_L,V) = \Deg(-A+F,V).
\ee
\indent For $t\in (0,1]$, define $\Upsilon_t:\o V\times [0,1]\to E$ by
$\Upsilon_t (x,\lma): = u(t; 0,1, x,\lma)$ where $u (t; 0,1, x,\lma)$
is a solution of $\dot u = - Au + \o F(u,\lma)$ on $[0,1]$ with the initial condition $u(0)=x$.
In view of Proposition \ref{semiflow-props-A-F}, $\Upsilon_t$ is well-defined and
$$
\beta \left( \Upsilon_t (\Omega\times [0,1]) \right) \leq e^{(k-\omega)t} \beta(\Omega)
\ \mbox{ for any } \Omega\subset \o V.
$$
In view of (\ref{degree-formula-III-step-4-2}) and Lemma \ref{degree-formula-II-lemma},
there exists $\widetilde t \in (0,1]$ such that, for any $t\in (0,\widetilde t]$,
\be\label{degree-formula-III-step-5-1}
\Upsilon_t (x,\lma) \neq x \ \mbox{ for any } (x,\lma)\in \part V\times [0,1].
\ee
Using the homotopy invariance of the degree, one has
\be\label{degree-formula-III-step-5-2}
\deg ( I - \Phi_t , V ) = \deg ( I - \Upsilon_t (\cdot,0) , V )
= \deg ( I - \Upsilon_t (\cdot,1) , V).
\ee
On the other hand, by Step 1, there exists $\o t \in (0, \widetilde t]$,
such that, for each $t\in (0,\o t]$,
$$
\deg (I-\Upsilon_t (\cdot,1), V) = \Deg ( -A+F_L , V),
$$
which together with (\ref{degree-formula-III-step-5-2}) and (\ref{degree-formula-III-step-4-3}) ends
the proof. \hfill $\square$


\begin{Rem} {\em
Theorem \ref{degree-formula-II} is a version of Th. 5.1 from \cite{Cwiszewski-1}, where $-A$ was assumed to generate a compact $C_0$ semigroup and $F$ to be locally Lipschitz with sublinear growth (and without any compactness properties). }
\end{Rem}

In the rest of the section we extend Theorem \ref{degree-formula-II} to the case
$k=\omega$, i.e. we  assume that (H1) and (H2) hold and, instead of (H3), we suppose that\\
\\ \noindent {\bf (H3a)}
\parbox[t]{120mm}{for any bounded $\Omega\subset E$ with $\beta(\Omega)>0$,
$$
\beta (F(\Omega\times [0,1])) < \omega \beta (\Omega).
$$}\\
\begin{Th}\label{condensing-th}
If $A:D(A)\to E$ and $F:E\to E$ satisfies $(H1)$, $(H2)$, $(H3a)$
and, for some open bounded $V\subset E$, $0\not\in (-A+F)(\part V\cap D(A))$,
then there exists $\o t>0$ such that, for any $t\in (0,\o t]$ and
$x\in \o V$, $\Phi_t (x)\neq x$, and
$$
\Deg(-A+F,V) = \deg(I-\Phi_t, V).
$$
\end{Th}
\begin{Lem} \label{niezw-lem}
Let $T_n:E\to E, \ n\ge 1$ be a sequence of bounded linear operators such that, for any $x\in E$,
$(T_n x)$ is a Cauchy sequence. Then, for any bounded set $\{x_n \}_{n\geq 1}\subset E$,
$$
\beta(\{T_n x_n\}_{n\ge 1})
\leq \left(\limsup_{n\to\infty}\|T_n\|\right)  \beta(\{x_n\}_{n\geq 1}).
$$
\end{Lem}
\noindent{\bf Proof}.
Observe first that, by the Banach-Steinhaus uniform boundedness theorem, the sequence $(\|T_n\|)$ is bounded.
Next take any $\varepsilon > 0$ and choose $y_{1},y_{2},\ldots ,y_{k}\in E$ making a
$\beta(\{x_n\}_{n\geq 1})+\varepsilon$ net for $\{x_n\}_{n\ge 1}$.
There exists an integer $N\geq 1$ such that, for any $l,n\geq N$ and $i=1,\dots,k$,
$\|T_n y_{i}-T_l y_{i}\|\leq \varepsilon$ and
$\|T_n\|\leq g +\varepsilon$ with $g:=\limsup_{m\to\infty}\|T_m\|$. Hence, for any $n\geq N$,
\begin{eqnarray*}
\|T_n x_n-T_{N} y_{i}\|\leq\|T_n x_n-T_n y_{i}\|+\|T_n y_{i}-T_N y_{i}\|
\leq \|T_n\| \|x_n - y_{i } \| + \eps\\
< (g+\eps)\left(\beta(\{x_m\}_{m\geq 1})+\varepsilon\right) + \eps
\end{eqnarray*}
where $i\in \{1,\ldots, k \}$ is chosen so that $\|x_n - y_{i } \| < \beta(\{x_m\}_{m\geq 1})+\varepsilon$.
This means that $\beta(\{ T_n x_n\}_{n\geq 1}) = \beta(\{ T_n x_n\}_{n\geq N}) \leq  (g+\eps)\left(\beta(\{x_m\}_{m\geq 1})+\varepsilon\right) + \eps$.
Since $\eps>0$ was arbitrary, we are done. \hfill $\square$

\noindent {\bf Proof of Theorem \ref{condensing-th}}.
We reduce the proof to the $k$-set contraction case (i.e. the situation of Theorem \ref{degree-formula-II}). To this end consider the parameterized family of differential problems
\be\label{reduction-equation-condensing}
\dot u = - A u + (1-\lma)F(u), \ \ \ \lma\in [0,1].
\ee
Define $\Theta_t:\o V\times [0,1]\to E$, for $t\in (0,1]$, by $\Theta_t (x,\lma):= u(t; 0,1, x,\lma)$ where
$u(\cdot;0,1, x,\lma)$ is a unique mild solution of $(\ref{reduction-equation-condensing})$ starting from $x$ at time $0$.
In view of Proposition \ref{semiflow-props-A-F}, $\Theta_t$ is well-defined, continuous
and
\be\label{condensing-th-1}
\beta( \Theta_t (\Omega\times [0,1])) \leq \beta(\Omega)
\mbox{ for any } \Omega\subset \o V.
\ee
\indent We shall now prove that,  for any $t\in (0,1]$,
$\Theta_t$ is a condensing map, i.e.
\be\label{condensing-th-2}
\beta \left( \Theta_t (\Omega\times [0,1]) \right) < \beta (\Omega)
\mbox{ for any }\Omega \subset \o V \mbox{ with } \beta(\Omega)>0.
\ee
Take any $\Omega$ with $\beta(\Omega)>0$.
Let $\Omega_0\subset \Omega$ be a countable set such that $\Omega \subset \o{\Omega_0}$
and $\Lambda_0 \subset [0,1]$ be countable and dense in $[0,1]$
(we use here the density argument as in the proof of Proposition \ref{semiflow-props-A-F}).
Then, clearly, for any $t\in (0, 1]$, by Proposition \ref{integral-prop-of-measure},
\begin{eqnarray*} 
\beta(\Theta_t(\Omega_0 \times \Lambda_0))  \ \ \ \ \ \mbox{ }
\ \ \ \ \ \ \ \ \ \ \ \ \ \ \ \ \ \ \ \ \ \ \ \ \ \ \ \ \ \ \ \ \ \ \ \ \ \ \ \ \ \ \ \ \ \ \ \ \ \ \ \ \ \ \
\mbox{ }\mbox{ }\mbox{ }\mbox{ }\mbox{ }\mbox{ }\mbox{ }\mbox{ }\mbox{ }\mbox{ }\mbox{ }\mbox{ }\mbox{ }\mbox{ }\mbox{ }\mbox{ }\mbox{ }\mbox{ }\mbox{ }  \\
\leq \nonumber
e^{-\omega t}\beta(\Omega_0) +
\beta\left(\left\{\left.\int_{0}^{t}  (1-\lma) S_A(t-\tau) F(\Theta_\tau(x,\lma)) \d\tau \right|
x\in \Omega_0, \, \lma \in  \Lambda_0\right\}\right)\\
\leq \nonumber
e^{-\omega t}\beta(\Omega_0) +
\beta\left(\left\{\left.\int_{0}^{t}  \gamma S_A(t-\tau) F(\Theta_\tau(x,\lma)) \d\tau \right|
x\in \Omega_0, \, \gamma,\lma \in  \Lambda_0\right\}\right)\\
\leq \nonumber e^{-\omega t}\beta(\Omega_0) +
\beta\left(\left\{\left.\int_{0}^{t}  S_A(t-\tau) F(\Theta_\tau(x,\lma)) \d\tau \right|
x\in \Omega_0, \, \lma \in  \Lambda_0\right\}\right)\\
\leq e^{-\omega t} \beta(\Omega_0) +
\int_{0}^{t}  \beta\left( S_A(t-\tau) F(\Theta_\tau(\Omega_0 \times \Lambda_0)) \right)\d\tau \\
\leq e^{-\omega t} \beta(\Omega_0) \nonumber + \int_{0}^{t} e^{-\omega(t-\tau)}
\beta(F(\Theta_\tau(\Omega_0 \times \Lambda_0))) \d\tau. \nonumber
\end{eqnarray*}
Further, either $\beta(\left\{ \Theta_\tau (x_n,\lma_n)\right\}_{n\geq 1})=0$, for a.e. $\tau\in [0,t]$,
and (\ref{condensing-th-2}) holds immediately or $\beta(\left\{ \Theta_\tau (x_n,\lma_n)\right\}_{n\geq 1})>0$
on a subset $J\subset [0,t]$ of positive measure.  In the latter case, by use of (\ref{condensing-th-1}),
for any $\tau\in J$,
$$
\beta (F(\Theta_\tau (\Omega_0\times \Lambda_0))) <  \omega\beta (\Theta_\tau (\Omega_0\times \Lambda_0))
\leq \omega\beta(\Omega_0)
$$
and, consequently,
\begin{eqnarray*}
\beta (\Theta_t(\Omega\times [0,1])) = \beta(\Theta_t(\Omega_0 \times \Lambda_0)) \leq e^{-\omega t} \beta(\Omega_0) \nonumber + \int_{0}^{t} e^{-\omega(t-\tau)}  \beta(F(\Theta_\tau(\Omega_0 \times \Lambda_0))) \d\tau\\
< e^{-\omega t} \beta(\Omega_0) + \int_{0}^{t} e^{-\omega(t-\tau)} \omega \beta(\Omega_0) d\tau\\
= e^{-\omega t} \beta(\Omega_0)  + (1-e^{-\omega t}) \beta(\Omega_0) = \beta (\Omega_0) = \beta(\Omega),
\end{eqnarray*}
i.e. (\ref{condensing-th-2}) holds.\\
\indent Next we shall show that there exist $\lma_0 \in (0,1]$ and $\o t\in (0, 1]$ such that, for any
$t\in (0, \o t]$,
\be\label{Krasnosel-II-3}
\Theta_t (x,\lma)\neq x \ \ \mbox{ for } (x,\lma)\in\part V\times [0,\lma_0].
\ee
If it were not so, then there would exist $(x_n)\subset \part V$, $\lma_n\to 0^+$ and
$t_n\to 0^+$ such that $\Theta_{t_n} (x_n,\lma_n)=x_n$ for each
$n\geq 1$. This would give $\Theta_{jt_n}(x_n, \lma_n)=x_n$ for
any integer $j>0$ such that $jt_n\leq 1$, and
\be\label{Krasnosel-II-3A}
\Theta_{1/2 + r_n} (x_n, \lma_n) =x_n
\ee
where $r_n:=([(1/2)/t_n]+1) t_n - 1/2 <t_n$, for $n\geq 1$. By the Duhamel formula
\be\label{Krasnosel-II-3B}
\Theta_{1/2+r_n}(x_n,\lambda_n) = S_A (r_n)\Theta_{1/2}(x_n,\lma_n) + I_n
\ee
where
$$
I_n:= \int_{0}^{r_n} S_A(r_n-s)F(\Theta_{1/2+s} (x_n,\lma_n)) \d s.
$$
Since $r_n\to 0^+$, we get $I_n\to 0$ as $n\to + \infty$.
Hence either $\{x_n\}_{n\geq 1}$ is relatively compact or,
by (\ref{Krasnosel-II-3A}), (\ref{Krasnosel-II-3B}), Lemma \ref{niezw-lem} and (\ref{condensing-th-2}),
\begin{eqnarray*}
\beta(\{ x_n\}_{n\geq 1}) & = \beta(\{\Theta_{1/2+r_n}(x_n,\lma_n)\}_{n\geq 1})
\leq \beta(\{ S_A (r_n)\Theta_{1/2}(x_n,\lma_n) \}_{n\geq 1}) \\
& \leq  \beta(\{ \Theta_{1/2}(x_n,\lma_n) \}_{n\geq 1}) < \beta (\{ x_n\}_{n\geq 1}),
\end{eqnarray*}
a contradiction proving that $\{x_n\}_{n\geq 1}$ is relatively compact.
Without loss of generality we may assume that $x_n\to x_0 \in \part V$.
If we put $u_n:= u(\cdot; 0,1, x_n, \lma_n)$, then, by Proposition \ref{semiflow-props-A-F},
$u_n\to u_0$ in $C([0, 1],E)$, and, using the $t_n$-periodicity of $u_n$,
we infer that $u_0\equiv x_0$, for some $x_0\in \part V\cap D(A)$,
and $-A x_0 + F(x_0)=0$, a contradiction completing the proof of (\ref{Krasnosel-II-3}).\\
\indent Finally, decreasing $\lma_0>0$ if necessary, we get
\be\label{Krasnosel-II-4}
\Deg (-A + F,V) = \Deg ( -A + (1-\lma_0)F , V )
\ee
and, in view of Theorem \ref{degree-formula-II} and  (\ref{Krasnosel-II-3}), there exists $\o t \in
(0, 1]$ such that, for any $t\in (0, \o t]$,
$
\Deg ( -A + (1-\lma_0)F , V ) = \deg( I-\Theta_t ( \cdot ,\lma_0) , V )
= \deg(I-\Theta_t ( \cdot ,0), V),$
which, together with (\ref{Krasnosel-II-4}), ends the proof. \hfill $\square$

\section{Periodic problem}

In this section we shall apply the degree formula from Theorem \ref{degree-formula-II} to the periodic problem of the form
$$
\left\{
\begin{array}{l}
\dot u = - Au + F(t,u) \mbox{ on } [0,T]\\
u(0)=u(T)
\end{array}
\right. \leqno{(P_T)}
$$
where
$A: D(A)\to E$ and $F:[0,T]\times E\to E $ satisfy conditions (H1)-(H2) and the following version of $(H3)$

\noindent {\bf (H3b)} \parbox[t]{120mm}{there exist $k\in
[0,\omega)$ such that, for any bounded $\Omega\subset E$,
$$
\beta(F([0,T]\times \Omega))\leq k\beta(\Omega).
$$}\\[1em]
Our approach is based on the idea of branching periodic solutions from the equilibrium point
of the the right-hand side (see \cite{Furi-Pera} and \cite{Cwiszewski-1})
in the parameterized differential equation. To this end consider the following parameterized family of equations
$$\dot u = -\lma Au + \lma F(t,u) \ \  \mbox{ on  }  [0, T],\ \lma\geq 0.
\leqno{(P_{T,\lma})}$$
By $u(\cdot;x,\lma)$ denote the unique mild solution of
$(P_{T,\lma})$ with the initial condition $u(0)=x$. Recall that a point
$(x,\lma)\in E \times [0,\infty)$ is called {\em a $T$-periodic point} for
$(P_{T,\lma})$ if $u(0;x,\lma)=u(T;x,\lma)$. Let $U\subset E$ be open and bounded.
One says that $x_0\in \o U$ is a {\em branching point} (or {\em cobifurcation point})
if there exists a sequence of $T$-periodic points $(x_n,\lma_n)\in \o U \times (0,+\infty)$
converging to $(x_0,0)$.
\begin{Th} \label{branching-ness-cond}
Let $A$, $F$ and $U$ be as above. If $x_0\in \o U$ is a branching point, then
$x_0\in D(A)$
and
$$
- A x_0 + \widehat F (x_0) = 0
$$
where $\widehat{F}: E \to E$ is given by
$\displaystyle{\widehat{F}(x):=\frac{1}{T}\int_{0}^{T} F(t,x)\d t}.$
\end{Th}
\begin{Rem} \label{rem-after-per-sol}{\em
Clearly $\widehat F$ is continuous and, combining Proposition \ref{integral-prop-of-measure} and assumptions (H2), (H3b), one has
$\beta(\widehat{F}(\Omega))\leq k\beta(\Omega) $ for any $\Omega\subset \o U.$ }
\end{Rem}
Since the Proof of Theorem \ref{branching-ness-cond} goes along the lines of
\cite[Th. 5.1]{Cwiszewski-1}, we skip it.

The Krasnosel'skii type degree formula allows us to derive a continuation principle.
\begin{Th} \label{cont-principle} {\em (Continuation Principle)}\\
Assume that $(H1), (H2)$ and $(H3b)$ hold and, for any $x\in E$,
$F(0,x)=F(T,x)$. If $(P_{T,\lma})$ has no periodic points in $\part U \times (0,1)$ and
$\Deg(-A+\widehat F, U)\neq 0$, then  $(P_T)$ admits a solution with $u(0)=u(T)\in \o U$.
\end{Th}
The proof of the theorem is based on the following averaging formula.
\begin{Prop} \label{mean-principle-prop}
Under the assumptions of Theorem \ref{cont-principle}, if $0\not\in(-A+\widehat{F})(\part U\cap D(A))$, then there exists
$\lma_0 > 0$ such that, for any $\lma \in (0 , \lma_0]$,
$$
\deg ( I - \Phi_T^\lma , U ) = \Deg (-A+\widehat{F},U) ,
$$
where $\Phi_T^\lma:\o U \to E$ is given by $\Phi_T^\lma(x):= u(T;x,\lma)$.
\end{Prop}
\noindent{\bf Proof.}   
Let $\widetilde{F}:[0,+\infty)\times E \to E$ be a continuous extension of $F$ given by $\widetilde{F} (t,x):=F(t-[t/T]T,x)$. And, for any $\lma\geq 0$,  define  $\Psi_t^\lma:\o U \times [0,1]\to E$ by
$\Psi_t^\lma(x,\mu):=u(t;x,\mu,\lma)$, for $(x,\mu)\in \o U \times [0,1]$, where
$u(\cdot;x,\mu,\lma)$ is a mild solution of
$$
\left\{
\begin{array}{l}
\dot u = -\lma A u + \lma \o F(t,u,\mu) \ \mbox{ on } [0,+\infty) \\
u(0)=x\\
\end{array}
\right.
$$
with $\o F: [0,+\infty)\times E\times [0,1]\to E$ given by
$\o F(t,x,\mu):=(1-\mu)\widehat{F}(x) + \mu \widetilde{F}(t,x).$\\
\indent Further note that, for any bounded $\Omega\subset E$,
\begin{eqnarray*}
\beta(\o F([0,+\infty)\times\Omega\times [0,1]))\leq
\beta(\cl{\conv}(\lma\widehat{F}(\Omega)\cup\lma {F}([0,T]\times\Omega)))=
\\ \lma \max\left\{\beta(\widehat{F}(\Omega)),\, \beta(\widetilde{F}([0,T]\times\Omega)\right\}
\leq \lma k \beta(\Omega),
\end{eqnarray*}
by (H2), for any $(t,x,\mu)\in [0,+\infty)\times E \times [0,1]$,
$$
\|\o F(t,x,\mu)\|\leq ((1-\mu)/T \int_0^T c \d t + \mu c)(1+\|x\|) = c(1+\|x\|).
$$
and $\|S_{\lma A}(t)\| \leq e^{-\lma\omega t }$ for any $t\geq 0$ and $\lma\geq 0$.
Hence, by Proposition \ref{semiflow-props-A-F} , $\Psi_{t}^{\lma}$ is continuous and
$$
\beta\left(\Psi_t^\lma (\Omega\times [0,1]) \right)\leq e^{\lma(k-\omega)t} \beta(\Omega) \mbox{ for any } \Omega \subset \o U.
$$
\indent We claim that there exist $\lma_1 > 0$ such that, for any $\lma\in (0,\lma_1]$,
\be\label{mean-principle-prop-0-1}
\Psi_T^\lma (x,\mu)\neq x, \mbox{ for any } (x,\mu) \in \part U \times [0,1].
\ee
If we suppose to the contrary, then there exist
$(\lma_n)\subset (0,+\infty)$, $(\mu_n) \subset [0,1]$ and $(x_n) \subset\part U$ such that $\lma_n\to 0^+$ and
\be\label{mean-principle-prop-0}
\Psi_T^{\lma_n} (x_n, \mu_n)=x_n  \mbox{ for each } n\geq 1.
\ee
Since, $\o{F}(\cdot, x,\mu)$ are $T$-periodic,
we get, for any $t>0$,
$$
\begin{array}{l}
\Psi_t^{\lma_n}(x_n,\mu_n) = \Psi_{t+T}^{\lma_n}(x_n, \mu_n) = \\
 = S_A (\lma_n T)  \Psi_t^{\lma_n} (x_n,\mu_n)
+\lma_n \displaystyle{\int\limits_t^{t+T}}
S_A(\lma_n (T+t-\tau)) \o F(\tau,\Psi_\tau^{\lma_n}(x_n,\mu_n),\mu_n)\d\tau
\end{array}
$$
and, consequently, for any integer $k\geq 1$,
\be\label{mean-principle-prop-1}
\begin{array}{l}
S_A((k-1)\lma_n T)\Psi_t^{\lma_n}(x_n,\mu_n)= S_A(k \lma_n T) \Psi_t^{\lma_n} (x_n,\mu_n)+\\
+ \lma_n S_A((k-1)\lma_n T) \displaystyle{\int\limits_t^{t+T}}
S_A (\lma_n(T+t-\tau)) \o F(\tau,\Psi_\tau^{\lma_n}(x_n,\mu_n),\mu_n)\d\tau.
\end{array}
\ee
If we put $k_n:=[1/\lma_n]+1$ and sum up the equalities (\ref{mean-principle-prop-1}) with $k=1,\ldots, k_n$, then we obtain
\be\label{mean-principle-prop-2}
\begin{array}{l}
  \Psi_t^{\lma_n} (x_n,\mu_n)  = S_A(k_n \lma_n T) \Psi_t^{\lma_n} (x_n,\mu_n)+\\
  \ \ \   +   \lma_n\left( \displaystyle{\sum\limits_{k=1}^{k_n}}
S_A(\lma_n (k-1) T)\right) \displaystyle{\int\limits_t^{t+T}}
S_A(\lma_n (T+t-\tau))\o F(\tau,\Psi_\tau^{\lma_n}(x_n,\mu_n),\mu_n)\d\tau  \\
\ \ \  \ \ \  =  S_A ((\lma_n k_n-1)T)S_A(T)\Psi_t^{\lma_n}(x_n,\mu_n) +
R_n \left( \displaystyle{\frac{1}{T} \int_{t}^{t+T} w_n(\tau) \d \tau} \right).
\end{array}
\ee
where $R_n := \lma_n T \displaystyle{\sum\limits_{k=1}^{k_n}}
S_A(\lma_n (k-1) T)):E\to E$
and $w_n:[t,t+T]\to E$, $n\geq 1$, are given by
$$
w_n(\tau):=S_A(\lma_n (T+t-\tau))\o F(\tau,\Psi_\tau^{\lma_n}(x_n,\mu_n),\mu_n).
$$
It is clear that, $\lma_n k_n \to 1$ and, for any $x\in E$, $R_n x\to \int\limits_{0}^{T} S_A (\tau) x \d \tau$ and
$$
\limsup_{n\to+\infty} \|R_n\| \leq
\limsup_{n\to+\infty} \sum_{k=1}^{k_n} e^{-\omega\lma_n(k-1)T } \lma_n T =
\int_{0}^{T} e^{-\omega \tau}\d \tau = \omega^{-1} (1-e^{-\omega T}).
$$
Hence, by use of (\ref{mean-principle-prop-2}), Proposition \ref{integral-prop-of-measure} and Lemma \ref{niezw-lem}, we have
$$
\begin{array}{l}
\beta\left(\{\Psi_t^{\lma_n}(x_n,\mu_n)\}_{n\geq 1}\right)
\displaystyle{\leq e^{-\omega T } \beta\left(\{\Psi_t^{\lma_n}(x_n,\mu_n)\}_{n\geq 1}\right)+ }\\
+ \displaystyle{\frac{1-e^{-\omega T}}{\omega T}}
\beta\left(\left\{\int_t^{t+T}w_n(\tau )\d\tau\right\}_{n\geq 1} \right)\\
\leq e^{-\omega T } \beta(\{\Psi_t^{\lma_n}(x_n,\mu_n)\}_{n\ge 1}) +
\displaystyle{\frac{1-e^{-\omega T} }{\omega T}} \int_t^{t+T}\beta(\{w_n (\tau )\}_{n\geq 1}) \d\tau \\
\leq e^{-\omega T } \beta(\{\Psi_t^{\lma_n}(x_n,\mu_n)\}_{n\geq 1}) +
\displaystyle{\frac{1-e^{-\omega T} }{\omega T}  \int_t^{t+T}
k  \beta (\{ \Psi_\tau^{\lma_n}(x_n,\mu_n) \}_{n\geq 1})\d\tau},
\end{array}
$$
which implies
\be\label{mean-principle-prop-3}
\beta(\{\Psi_t^{\lma_n}(x_n,\mu_n)\}_{n\ge 1})\leq
\displaystyle{\frac{k}{\omega T} \int_t^{t+T}\beta(\{\Psi_\tau^{\lma_n}(x_n,\mu_n)\}_{n\geq 1})\d\tau}.
\ee
Now define $\varphi:[0,+\infty)\to \R$ by
$\varphi(\tau):=\beta(\{\Psi_\tau^{\lma_n}(x_n,\mu_n)\}_{n\geq 1})$. By Remark \ref{semiflow-gronwall},
it is clear that $\varphi$ is bounded and, consequently, $M:=\sup_{\tau \in [0,T]} \varphi(\tau) <+\infty$.
If $M>0$, then there exists $0<\eps<M(1-k/\omega)$ and $t_\eps\geq 0$ such that $\varphi(t_\eps) > M-\eps$.
Then, (\ref{mean-principle-prop-3}), gives
$$
M-\eps < \varphi(t_\eps) \leq \displaystyle{\frac{k}{\omega T} \int_{t_\eps}^{t_\eps+T}\varphi(\tau)\d\tau} \leq (k/\omega) M < M-\eps,
$$
a contradiction meaning that $M=0$. Hence $\beta(\{x_n\}_{n\ge 1}) = 0$.
Hence without loss of generality we may assume that $x_n\to x_0$ and $\mu_n\to \mu_0$
for some $x_0\in\part U$ and $\mu_0\in [0,1]$. If we put
$u_n:=u(\cdot, x_n, \mu_n, \lma_n )$, then, obviously, for any $t\in [0,T]$,
$$
u_n(t)= S_A(\lma_n t)x_n + \lma_n\int_0^t S_A(\lma_n(t-\tau)) \o F(\tau, u_n(\tau),\mu_n)\d \tau
$$
and, further,
$$
\|u_n(t)-x_0\|\leq \|S_A(\lma_n t) x_n - x_0\| + \lma_n c T(1+K)
$$
where $K$ is the bound for solutions starting from $\o U$.
Passing to the limit, one sees that $(u_n)$ converges
uniformly on $[0,T]$ to $u_0 \equiv x_0$.
Therefore, using again (\ref{mean-principle-prop-0}), we obtain
$$
-\frac{S_A(\lma_n T)x_n -x_n}{\lma_n T} =
\frac{1}{T} \int_0^T S_A(\lma_n(T-\tau)) \o F(\tau,u_n(\tau),\mu_n) \d\tau,
$$
and, in consequence,
$$
A\left( \frac{1}{\lma_n T } \int_{0}^{\lma_n T} S_A(\tau)x_n \d \tau \right) = \frac{1}{T} \int_0^T S_A(\lma_n(T - \tau)) \o F(\tau,u_n(\tau),\mu_n) \d\tau.
$$
Since $\frac{1}{\lma_n T}\int_0^{\lma_n T}S_A(\tau) x_n \d \tau\to x_0$,
$$
\frac{1}{T}\int_0^T S_A(\lma_n(T-\tau))\o F(\tau,u_n(\tau),\mu_n) \d\tau \to
\frac{1}{T}\int_0^T \o F(\tau,x_0,\mu_0)\d \tau = \widehat{F}(x_0),
$$
and $A$ is closed, we get  $-Ax_0+\widehat{F}(x_0) = 0$, which is a contradiction ending the proof of (\ref{mean-principle-prop-0-1}).\\
\indent Thus, using the homotopy invariance of the topological degree, we get, for any
$\lma\in (0,\lma_1]$,
\be\label{mean-principle-prop-4}
\begin{array}{ll}
\deg(I -\Phi_T^\lma , U ) & = \deg (I-\Psi_T^{\lma}(\cdot,1),U)=\deg(I -\Psi_T^\lma(\cdot,0),U)\\
& =\deg(I-\widehat{\Phi}_{T}^{\lma},U)
\end{array}
\ee
where $\widehat{\Phi}_{t}^{\lma}$ is the translation along trajectory operator for the equation
$\dot u = -\lma A u +\lma\widehat{F} (u)$. In view of Theorem \ref{degree-formula-II},
there exists $\lma_0\in (0,\lma_1]$ such that, for any $\lma\in (0,\lma_0]$,
$\deg(I - \widehat{\Phi}_{\lma T}^{1} ,U) = \Deg(-A+\widehat{F},U),$
which, by rescaling the time, gives
\be\label{mean-principle-prop-5}
\deg(I-\widehat \Phi_{T}^{\lma}, U)=\deg(I - \widehat {\Phi}_{\lma T}^{1} ,U) = \Deg(-A+\widehat{F},U).
\ee
Combining (\ref{mean-principle-prop-4}) and (\ref{mean-principle-prop-5}) together,
we complete the proof. \hfill $\square$

\noindent {\bf Proof of Theorem \ref{cont-principle}}.
In view of Proposition \ref{mean-principle-prop}, there
exists $\lma_0>0$ such that, for $\lma\in (0,\lma_0]$,
$$
\deg (I-\Phi_{T}^{\lma}, U) = \Deg(-A+\widehat F, U).
$$
Either $\Phi_{T}^{1}$ has a fixed point in $\part U$ or,
by the assumption and the homotopy invariance of the degree,
$$
\deg(I-\Phi_{T}^{1},U) = \deg(I-\Phi_{T}^{\lma_0},U)=\Deg(-A+\widehat{F},U)\neq 0.
$$
and there exists $x\in U$ such that $(x,1)$ is a periodic point.
In both cases one gets the exists of a $T$-periodic solution
starting (and ending) at a point of $\o U$.     \hfill $\square$

\begin{Rem}{\em
(a)  If $U:=B(x_0,r)$ for some $x_0\in D(A)$ and $r>0$ such that  $-A x_0 + \widehat F(x_0) = 0$, and $\Deg(-A+\widehat F, B(x_0,r))\neq 0$ and $(P_{T,\lma})$ has no periodic points in $\part B(x_0,r) \times (0,1)$, then the thesis of Theorem \ref{cont-principle} follows also from Theorem 6.2.1 in \cite{ObuKaZ}. It is worth mentioning that in \cite{ObuKaZ} the authors work with set-valued $F$ and a periodic solution comes up as a fixed point of a certain operator in the space of periodic functions from $C([0,T],E)$.\\
\indent (b) For the set-valued upper semicontinuous $F$ with compact convex values,
the translation along trajectories approach also applies.  Obviously, in general, the inclusion $\dot u\in - A u + F(u)$ has not a unique solution, i.e. the translation along trajectories operator $\Phi_t:E\multimap E$ has the form
$$
\Phi_{t}^{(F)}(x)  = e_t (L(x,-A+F))
$$
where $L(x,-A+F)$ is the set of solution starting from $x$ and  $e_t: C([0,T],E)\to E$ is the evaluation at $t$. However, due to results of \cite{BaKr2}, the solution set $L(x,-A+F)$, for any $x\in E$, is cell-like (or of $R_\delta$ type).
This enables to apply the fixed point index theory for a general class of maps given by a superposition $e\circ L : E \to E$ where $L:E\to E_1$ is  a set-valued upper semicontinuous map with $R_\delta$ values and $e:E_1\to E$ is a continuous map -- see \cite{Kryszewski} and \cite{BaKr1}.
Using this degree and the related approximation method, we reduce the problem to the single-valued one.
Roughly speaking, taking sufficiently close graph approximation $f$ of $F$, we know that the topological degree of $-A+f$ is equal to that of $-A+F$ and
the degree of $I-\Phi_{t}^{(f)}$ is equal to the degree of $I-\Phi_{t}^{(F)}$.
}\end{Rem}

Applying Theorem \ref{cont-principle}, we obtain the following simple existence criterion.
\begin{Th}\label{comp-periodic-principle}
Suppose  $(H1)$ and $(H3b)$ hold, $F(0,x)=F(T,x)$ for any $x\in E$, and there exists $K>0$ such that
$$
\|F(t,x)\|\leq K \mbox{ for any } (t,x)\in [0,T]\times E,
$$
then the periodic problem $(P_T)$ admits a mild solution $u$ with $\|u(0)\|=\|u(T)\|\leq K / \omega$.
\end{Th}
\noindent {\bf Proof}. We claim that $(P_{T,\lma})$ has no periodic points in
$\part B(0,R)\times (0,1)$ if $R> K / \omega$.
Indeed, suppose to the contrary, that there exists a periodic point
$(\o x,\o \lma)\in \part B(0,R)\times (0,1)$.
Then
$$
\o x = S_{\o\lma A }(T)\o x + \o\lma \int_{0}^{T} S_{\o\lma A }(T-\tau)
F(\tau, \o u(\tau)) \d\tau
$$
where $\o u$ is the corresponding periodic mild solution of $(P_{T,\o \lma})$.
By the assumption, we get
$$
\| \o x \| \leq e^{-\o\lma\omega T} \|\o x\| + \o \lma K \int_{0}^{T}
e^{-\o\lma\omega (T-\tau)} \d\tau,
$$
which gives
$R\leq Re^{-\o\lma\omega T} (1+\o\lma K/R (\o\lma \omega)^{-1}(e^{\o\lma\omega T} - 1))
= R e^{-\o\lma\omega T} (1+ K/(R\omega)(e^{\o\lma\omega T} - 1))<R$,
a contradiction.\\
\indent To compute the degree $\Deg(-A+\widehat F,B(0,R))$ consider
the homotopy $A  + \mu \widehat F$. It is admissible, since
$\beta( \bigcup_{\mu\in [0,1]} \mu \widehat F (\Omega))\leq k \beta(\Omega)$ for any bounded $\Omega\subset E$.
If there exists $(x,\mu)\in \part B(0,R)\cap D(A)\times [0,1]$ such that
$-Ax+\mu \widehat F(x)=0$, then $x = \mu A^{-1} \widehat F(x)$
and, in consequence,
$R=\|x\|\leq \mu \|A^{-1}\| K \leq \mu K/\omega <  R$, a contradiction.
Hence, by the homotopy invariance of the topological degree
$$
\Deg(-A+\widehat F, B(0,R)) = \Deg(-A,B(0,R)) = \deg (I,B(0,R))=1
$$
Applying Theorem \ref{cont-principle}, we get the desired periodic solution in $\o{B(0,R)}$.\\
\indent Since the translation along trajectory operator $\Phi_T$ is a $k$-set contraction and
$R$ was arbitrary, we can conclude that $(P_T)$ has a solution in $\o{B(0,K/\omega)}$. \hfill $\square$

We end the paper with an example of applications to a system of partial differential equations.
\begin{Ex}{\em
Consider the following system of equations (cf. \cite{Hale} and \cite{ObuKaZ})
coming from the theory of transmission lines
\be\label{example-0}
\, \left\{
\begin{array}{l}
\displaystyle{\frac{\part I}{\part t}(t,x) + \alpha\displaystyle{\frac{\part V}{\part x}}(t,x)+ \beta I(t,x) =f(t,x,I(t,x),V(t,x)) \  x\in [0,l], t\geq 0 }\\
\displaystyle{\frac{\part V}{\part t}}(t,x) +\gamma \displaystyle{\frac{\part I}{\part x}}(t,x)
+ \delta V(t,x) = g(t,x,I(t,x),V(t,x))  \  x\in [0,l], t\geq 0\\
V(t,0)+\rho I(t,0) = 0 \ \ \ \ \ \ \ \   t\geq 0\\
-I(t,l)+\sigma \displaystyle{\frac{\part V}{\part t}}(t,l)  + h(V(t,l)+e(t))+j(t)=0 \ \ \ \ \ t\geq 0
\end{array}
\right.
\ee
where $\alpha,\beta,\gamma,\delta,\rho,\sigma>0$, continuous $j, e :\R\to\R$ are $T$-periodic,  $f,g:\R\times [0,l]\times \R^2\to\R$ are continuous, bounded, $T$-periodic with respect to the first variable and there are $L_f, L_g\geq 0$ such that, for any $(t,x)\in \R\times [0,l]$ and $(v_1, v_2), (w_1, w_2)\in\R^2$,
$$
\begin{array}{ll}
|f(t,x,v_1,w_1)-f(t,x,v_2,w_2)|\leq L_f|(v_1,w_1)-(v_2,w_2)| & \mbox{ and }\\
|g(t,x,v_1,w_1)-g(t,x,v_2,w_2)|\leq L_g|(v_1,w_1)-(v_2,w_2)| &
\end{array}
$$
($|\cdot|$ stands for the Euclidean norm in $\R^2$), $h:\R\to\R$  is Lipschitz with a constant $L_h\geq 0$
and with the property that there exist $a>\sigma \min\{ \beta,\delta \},\, b\geq 0$ such that
\be\label{example-1}
| h (s) - a s | \leq b \ \ \mbox{ for any } s \in \R
\ee
and
\be \label{example-2}
\max \left\{ L_f^2+\alpha L_g^2/\gamma, \,\gamma L_f^2/ \alpha +L_g^2,\,  2(a^2+L_h^2)/\sigma \right\} < (\min\{\beta,\delta\})^2.
\ee
Let  $E:=L^2 (0,l)\times L^2(0,l)\times \R$ be a linear space equipped with the scalar product given by
$$\la (p_1,q_1,r_1), (p_2,q_2, r_2)\ra_E:= \gamma\la p_1, p_2 \ra+\alpha\la q_1, q_2\ra +
\alpha \gamma \sigma r_1 r_2$$
where $\la \cdot,\cdot \ra$ stands for the scalar product in $L^2(0,l)$
Define $A:D(A)\to E$ by
$$A(p,q,r):= \left(
\begin{array}{lll}
\beta & \alpha \part_x & 0\\
\gamma \part_x & \delta & 0 \\
-\sigma^{-1} e_l & 0 & \sigma^{-1}a
\end{array}
\right)
\left( \begin{array}{l}
p\\
q\\
r\\\end{array}\right),   (p,q,r)\in D(A)
$$
where
$D(A):= \{ (p,q,r)\in E\,|\, p,q \in H^1(0,l),\, q(0)+\rho p(0)=0,\, q(l)=r  \}$, $e_l$ is the evaluation
at $l$ operator on $H^1(0,l)$. Observe that, for any $(p,q,r)\in D(A)$, integrating by parts
$$
\begin{array}{l}
\la A(p,q,r), (p,q,r) \ra_E =
\gamma \la \beta p + \alpha q', p\ra + \alpha \la \gamma p'+\delta q, q\ra + \alpha \gamma \sigma (-\sigma^{-1})(p(l)-a r)r  \\
=\beta \gamma \|p\|_{L^2}^{2} +
\alpha \gamma (qp|_{0}^{l}-\la p',q\ra) + \alpha\delta\|q\|_{L^2}^{2} +
\alpha \gamma \la p',q\ra -\alpha\gamma p(l)r + \alpha \gamma a r^2\\
= \beta \gamma \|p\|_{L^2}^{2} + \alpha\delta\|q\|_{L^2}^{2} + \alpha \gamma a r^2
- \alpha \gamma q(0)p(0) \geq \min\{\beta, \delta, a/\sigma\}
\|(p,q,r)\|_{E}^{2} + \alpha\gamma \rho p(0)^{2}\\
\geq  \min\{\beta, \delta, a/\sigma \} \|(p,q,r)\|_{E}^{2}.
\end{array}
$$
Moreover, by direct computation we can see that the range of $A+\lma I$, for $\lma>0$, is equal to $E$.
Hence, due to the Lumer-Philips theorem, $-A$ generates a $C_0$ semigroup $S_A$ on $E$ and   $\| S_A(t)\|\leq \exp (- \min\{\beta, \delta, a/\sigma\} t)$
for any $t\geq 0$, i.e. $(H1)$ is satisfied.\\
\indent If we define $F:[0,+\infty)\times E \to E$ by
$$
F(t,p,q,r):= \left(
\begin{array}{c}
N_f(t,p,q)\\ 
N_g(t,p,q)\\ 
\sigma^{-1}(ar- h(r+e(t))-j (t)),
\end{array}
\right).
$$
where $N_f, N_g$ are given by $N_f(t,p,q)(x):=f(t,x,p(x),q(x))$, $N_g(t,p,q)(x):=f(t,x,p(x),q(x))$, for
a.e. $x\in [0,l]$, then it is clear that (\ref{example-0}) can be put as an abstract equation (with $u:=(p,q,r)$)
$$
\dot u = - Au +F(t,u).
$$
We shall apply Theorem \ref{comp-periodic-principle} to find a periodic (mild) solution of the system (\ref{example-0}). To this end, observe that, under these assumptions, for any $(p_1,q_1, r_1), (p_2,q_2, r_2)\in E$ and $t>0$,
$$
\begin{array}{l}
\|F(t,p_1,q_1,r_1) - F(t,p_2,q_2, r_2)\|_{E}^{2}
= \gamma \|N_f(t,p_1,q_1)- N_f(t,p_2,q_2)\|_{L^2}^{2} +\\
+ \alpha \| N_g(t,p_1,q_1)- N_g(t,p_2,q_2)\|_{L^2}^2+ \alpha \gamma \sigma \left|
\sigma^{-1}(a(r_1-r_2) - (h(r_1+e(t))-h(r_2+e(t))))\right|^2 \\
\leq (\gamma L_f^2 + \alpha L_g^2)(\|p_1-p_2\|_{L^2}^{2} + \|q_1-q_2\|_{L^2}^{2})
  + \alpha \gamma  2(a^2+L_h^2)|r_1-r_2|^{2} \\
  \leq \max
\left\{ L_f^2+\alpha L_g^2/\gamma, \,\gamma L_f^2/ \alpha +L_g^2,\,
  2(a^2+L_h^2)/\sigma \right\}\|(p_1,q_1,r_1)-(p_2,q_2,r_2)\|_{E}^{2},
\end{array}
$$
which, in view of (\ref{example-2}), means that $(H3b)$ is satisfied.
It is also easy to see that $F$ is bounded, since first two coordinates are bounded directly by the assumption
and, for the third one,
\begin{eqnarray*}
|\sigma^{-1}(ar- h(r+e(t))-j (t))| \leq \sigma^{-1}(|a(r+e(t))-h(r+e(t))| + a|e(t)| + |j(t)|) \\
\leq
\sigma^{-1} (b+a\max_{\tau\in [0,T]} |e(\tau)| + \max_{\tau\in [0,T]} |j(\tau)|).
\end{eqnarray*}
Thus, by use of Theorem \ref{comp-periodic-principle}, we conclude that (\ref{example-0}) possesses a $T$-periodic (mild) solution.
}
\end{Ex}
We finish with an example showing that the results of the paper may be also applied to
parabolic partial differential equations where the semigroup of the proper differential operator
may not be compact and the results of \cite{Cwiszewski-1} do not apply.
\begin{Ex}{\em
Suppose $\Omega\subset \R^N$ is an open set contained in a strip of with $d>0$, i.e.
there is $j_0\in \{1,\ldots, N\}$ such that, for any $x\in \Omega$, $|x_{j_0}|\leq d$.
Define $A:D(A)\to L^2(\Omega)$ by $D(A):=\{ u\in H_{0}^{1}(\Omega)\,|\, \Delta u \in L^2(\Omega) \}$
and $Au := -\Delta u$. Then $-A$ generates a $C_0$ semigroup $S_A$ on $L^2(\Omega)$ such that
$\|S_A (t)\|\leq e^{-\omega t}$ for some $\omega>0$ and all $t\geq 0$. Indeed,
by the Poincar\'{e} inequality,
$$
\la Au, u\ra_{L^2} = \int_{\Omega} |\nabla u|^{2} \d x \geq 4d^{-2} \|u\|_{L^2}^{2},
$$
which guarantees that $-A$ generates a $C_0$ semigroup such that $\|S_A(t)\|
\leq e^{-\omega t}$ with $\omega:=4 d^{-2}$.
}
\end{Ex}

\noindent {\bf Acknowledgments}
The authors would like to thank the autonomous referees for suggesting an improvement of Theorem
\ref{condensing-th} and indicating some important references.


\begin{thebibliography}{99}
\newcommand{\bi}{\bibitem}
\baselineskip 2 mm

\bi{Akhmerov-etal} Akhmerov R.R., Kamenskii M.I., Potapov A. S., Rodkina A.E., Sadovskii B.N.,
{\em Measures of noncompactness and condensing operators}, Birkh\"{a}user 1992.


\bi{BaKr1} Bader R., Kryszewski W., {\em Fixed-point index for compositions of set-valued maps with proximally $\infty$-connected values on arbitrary ANR's},  Set-Valued Anal.  2  (1994),  no. 3, 459--480.

\bi{BaKr2} Bader R.,  Kryszewski W., {\em On the solutions of differential inclusions
and the periodic problem in Banach spaces}, Nonlinear An. 54 (2003), 707--754.

\bi{Bothe}
 Bothe D. {\em Periodic solutions of a nonlinear evolution problem
 from heterogeneous catalysis}, Diff. Int. Eq.,  vol. 16, no. 6 (2001), 641--670.

\bi{Cwiszewski-1} \'{C}wiszewski A., {\em Topological degree methods for perturbations
                                     of operators generating compact $C_0$ semigroups},
                                     J. Diff. Eq., Vol. 220, No. 2 (2006), 434--477.
\bi{Dugundji-Granas} Dugundji J., Granas A., {\em Fixed Point Theory}, Springer, Berlin 2004.

\bi{Deimling} Deimling K., {\em Multivalued Differential Equations}, Walter de Gruyter, Berlin 1992.


\bi{Furi-Pera} Furi M., Pera M. P., {\em Global Bifurcation of Fixed
Points and the Poincar\'{e} Translation Operator on Manifolds},
Annali di Matematica pura ed applicata (IV), Vol. CLXXIII (1997),
pp. 313--331.


\bi{Hale} Hale, J. K., {\em Asymptotic behavior of dissipative systems},
Mathematical Surveys and Monographs 25, American Mathematical Society, Providence, RI, 1988.

\bi{ObuKaZ} Kamenskii M., Obukhovskii V., Zecca P., {\em Condensing Multivalued Maps and Semilinear
             Differential Inclusions in Banach Spaces},
             De Gruyter Series in Nonlinear Analysis and Applications 7, Walter de Gruyter, 2001.
\bi{Kartsatos} Kartsatos A. G. , {\em Recent results involving compact perturbations and
compact resolvents of accretive operators in Banach spaces},
                Proceedings of the First World Congress of Nonlinear Analysts,
                Tampa, Florida, 1992, Walter de Gruyter, New York 1995.

\bi{Krasnoselskii} Krasnosel'skii M. A., {\em The Operator of Translation along Trajectories of Ordinary
Differential Equations}, Translations of Mathematical Monographs, Vol. 19, AMS, Providence, R.I. 1968.

\bi{Kryszewski} Kryszewski W.,  {\em Graph-approximation of set-valued maps on noncompact domains},
 Topology Appl.  83 (1998), no. 1, 1--21.


\bi{Mawhin69} Mawhin, J., {\em \'{E}quations int\'{e}grales et solutions p\'{e}riodiques des syst\`{e}mes
differentiels non lin\'{e}ares}, Acad. Roy. Belg. Bull. Cl. Sci. (5) 55 1969, 934--947.

\bi{Nussbaum} Nussbaum, R. D., {\em The fixed point index for local condensing maps},
Ann. Mat. Pura Appl. (4) 89 (1971), 217--258.

\bi{Pazy} Pazy A., {\em Semigroups of linear operators and applications to partial
                        differential equations}, Springer Verlag 1983.

\bi{Schi-Schm} Schiaffino A., Schmitt K., {\em Periodic Solutions of a Class of Nonlinear Evolution Equations}, Ann. Mat. Pura Appl. (4) 137 (1984), 265--272.

\end{thebibliography}
\end{document}